\documentclass[10pt]{amsart}
\usepackage{amsmath, amssymb, amsthm}
\usepackage{amscd}
\usepackage[dvips]{graphicx}

\usepackage[all]{xy}
\usepackage[pdftitle={Orthogonal decomposition and Lehmer's problem},pdfauthor={Fili and Miner},colorlinks=false,pdfstartview={}]{hyperref}

\newtheorem{thm}{Theorem}

\newtheorem*{thm2}{Theorem \ref{thm:Lp-Lehmer-equiv}}
\newtheorem*{thm3}{Theorem \ref{thm:Gamma-equiv}}
\newtheorem*{thm4}{Theorem \ref{thm:mahler-2-norm}}
\newtheorem{conj}{Conjecture}
\newtheorem{prop}{Proposition}
\newtheorem{lemma}[prop]{Lemma}

\newtheorem{cor}[prop]{Corollary}

\newtheorem{example}[prop]{Example}

\newtheorem*{thm*}{Theorem}
\newtheorem*{alg*}{Algorithm}
\newtheorem*{lemma*}{Lemma}

\theoremstyle{remark}
\newtheorem{rmk}[prop]{Remark}
\newtheorem*{rmk*}{Remark}
\newtheorem{notation}[prop]{Notation}
\newtheorem*{notation*}{Notation}

\theoremstyle{definition}
\newtheorem{defn}[prop]{Definition}

\numberwithin{prop}{section}
\numberwithin{equation}{section}

\newcommand{\mybf}{\mathbb}

\newcommand{\bR}{\mybf{R}}

\newcommand{\bC}{\mybf{C}}
\newcommand{\bN}{\mybf{N}}
\newcommand{\bZ}{\mybf{Z}}
\newcommand{\bQ}{\mybf{Q}}

\newcommand{\cP}{\mathcal{P}}
\newcommand{\cB}{\mathcal{B}}

\newcommand{\cL}{\mathcal{L}}

\newcommand{\cF}{\mathcal{F}}

\newcommand{\cU}{\mathcal{U}}
\newcommand{\cG}{\mathcal{G}}
\newcommand{\cI}{\mathcal{I}}
\newcommand{\cK}{\mathcal{K}}

\newcommand{\al}{\alpha}

\newcommand{\Gal}{\operatorname{Gal}}

\newcommand{\supp}{\operatorname{supp}}

\newcommand{\ra}{\rightarrow}

\newcommand{\lcm}{\operatorname{lcm}}

\newcommand{\ep}{\epsilon}

\newcommand{\lfrac}[2]{\left(\frac{#1}{#2}\right)}

\newcommand{\Tor}{\operatorname{Tor}}

\newcommand{\Stab}{\operatorname{Stab}}
\newcommand{\ord}{\operatorname{ord}}

\newcommand{\Qbar}{\overline{\mybf{Q}}}

\def\talltareesidedbox#1{\setbox0=\hbox{$#1$}\dimen0=\wd0 \advance\dimen0 by3pt\rlap{\hbox{\vrule height10pt width.4pt
 depth2pt \kern-.4pt\vrule height10.4pt width\dimen0 depth-10pt\kern-.4pt \vrule height10pt width.4pt depth2pt}}
 \relax \hbox to\dimen0{\hss$#1$\hss}}
\def\tareesidedbox#1{\setbox0=\hbox{$#1$}\dimen0=\wd0 \advance\dimen0 by3pt\rlap{\hbox{\vrule height8pt width.4pt
 depth2pt \kern-.4pt\vrule height8.4pt width\dimen0 depth-8pt\kern-.4pt \vrule height8pt width.4pt depth2pt}}
\relax \hbox to\dimen0{\hss$#1$\hss}}
\def\shorttareesidedbox#1{\setbox0=\hbox{$#1$}\dimen0=\wd0 \advance\dimen0 by3pt\rlap{\hbox{\vrule height7pt width.4pt
 depth2pt \kern-.4pt\vrule height7.4pt width\dimen0 depth-7pt\kern-.4pt \vrule height7pt width.4pt depth2pt}}
 \relax \hbox to\dimen0{\hss$#1$\hss}}
\newcommand{\house}[1]{\tareesidedbox{#1}}
\newcommand{\tallhouse}[1]{\talltareesidedbox{#1}}
\newcommand{\shorthouse}[1]{\shorttareesidedbox{#1}}

\newcommand{\bal}{{\alpha}}

\newcommand{\fal}{f_\alpha}
\newcommand{\ovrln}{}

\newcommand{\N}{\operatorname{N}}

\title[Orthogonal decomposition and Lehmer's problem]{Orthogonal decomposition of the space of algebraic numbers and Lehmer's problem}
\author[Fili]{Paul Fili}
\address{Department of Mathematics\\ University of Texas at Austin, TX 78712}
\email{pfili@math.utexas.edu}
\author[Miner]{Zachary Miner}
\address{Department of Mathematics\\ University of Texas at Austin, TX 78712}
\email{zminer@math.utexas.edu}
\date{\today}

\begin{document}

\begin{abstract}
We introduce vector space norms associated to the Mahler measure by using the $L^p$ norm versions of the Weil height recently introduced by Allcock and Vaaler. In order to do this, we determine orthogonal decompositions of the space of algebraic numbers modulo torsion by Galois field and degree. We formulate $L^p$ Lehmer conjectures involving lower bounds on these norms and prove that these new conjectures are equivalent to their classical counterparts, specifically, the classical Lehmer conjecture in the $p=1$ case and the Schinzel-Zassenhaus conjecture in the $p=\infty$ case.
\end{abstract}

\maketitle

\tableofcontents

\section{Introduction}

Let $K$ be a number field with set of places $M_K$. For each $v\in M_K$ lying over a rational prime $p$, let $\|\cdot\|_v$ be the absolute value on $K$ extending the usual $p$-adic absolute value on $\bQ$ if $v$ is finite or the usual archimedean absolute value if $v$ is infinite. Then for $\al\in K^\times$, the \emph{absolute logarithmic Weil height} $h$ is given by 
\[
 h(\al) = \sum_{v\in M_K} \frac{[K_v:\bQ_v]}{[K:\bQ]} \log^+ \|\al\|_v 
\]
where $\log^+ t = \max\{0,\log t\}$. As the expression on the right hand side of this equation does not depend on the choice of field $K$ containing $\al$, $h$ is a well-defined function mapping  $\Qbar^\times\ra[0,\infty)$ which vanishes precisely on the roots of unity $\Tor(\Qbar^\times)$. Closely related to the Weil height is the \emph{logarithmic Mahler measure}, given by
\[
 m(\al)=(\deg_{\bQ} \al)\cdot h(\al)
\]
where $\deg_{\bQ}\al = [\bQ(\al):\bQ]$. Though seemingly related to the Weil height in a simple fashion, the Mahler measure is in fact a fair bit more mysterious. Perhaps the most important open question regarding the Mahler measure is \emph{Lehmer's problem}, which asks if there exists an absolute constant $c$ such that
\begin{equation}
 m(\al)\geq c>0\quad\text{for all}\quad \al\in\Qbar^\times\setminus \Tor(\Qbar^\times).
\end{equation}
The question of the existence of algebraic numbers with small Mahler measure was first posed in 1933 by D.H. Lehmer \cite{Lehmer} and since then the conjectured existence of an absolute lower bound away from zero has come to be known as \emph{Lehmer's conjecture}. The current best known lower bound, due to Dobrowolski \cite{Dob}, is of the form
\[
m(\al)\gg \lfrac{\log \log \deg_{\bQ} \al}{\log \deg_{\bQ} \al}^3\quad\text{for all}\quad \al\in\Qbar^\times\setminus \Tor(\Qbar^\times)
\]
where the implied constant is absolute. 

Recently, Allcock and Vaaler \cite{AV} observed that the absolute logarithmic Weil height $h:\Qbar^\times\ra [0,\infty)$ can in fact be viewed in an equivalent fashion as the $L^1$ norm on a certain measure space $(Y,\lambda)$. The points of $Y$ are the places of $\Qbar$ endowed with a topology which makes $Y$ a totally disconnected locally compact Hausdorff space, and each equivalence class of the algebraic numbers modulo torsion give rise to a unique locally constant function with compact support. The purpose of this paper is to construct analogous function space norms in order to study the Mahler measure. Once we have introduced our new norms, we will give $L^p$ forms of the Lehmer conjecture which are equivalent to the classical Lehmer conjecture for $p=1$ and to the Schinzel-Zassenhaus conjecture for $p=\infty$.

We first briefly recall here the notation of \cite{AV}, which we will use throughout this paper. To each equivalence class $\al$ in $\Qbar^\times/\Tor(\Qbar^\times)$, we can uniquely associate the function $\fal:Y\ra \bR$ given by \[\fal(y) = \log \|\al\|_y\] (we will often drop the subscript $\al$ when convenient). We denote the space of functions given by algebraic numbers modulo torsion by $\cF$. If $\al\in K$, then the function $\fal(y)$ is constant on the sets $Y(K,v)=\{y \in Y : y|v\}$ for $v\in M_K$ and takes the value $\log \|\al\|_v$. Then if $\al\in K^\times$ for some number field $K$, we have
\[
 \|\fal\|_1=\int_Y |\fal(y)|\,d\lambda(y) = \sum_{v\in M_K} |\log \|\al\|_v|\,\frac{[K_v:\bQ_v]}{[K:\bQ]} = 2\,h(\al).
\]
The product formula takes the form
$
 \int_Y \fal \,d\lambda=0.
$
We also have a well-defined inner product on $\cF$ given by
\[
 \langle f,g\rangle =\int_Y f(y)g(y)\,d\lambda(y)
\]
which satisfies $\|f\|_2 = \langle f,f\rangle^{1/2}$. The geometry of the space of $\cF$ will play a significant role in our study.

The study of the Mahler measure on the vector space of algebraic numbers modulo torsion $\cF$ presents several difficulties absent for the Weil height, first of which is that $m$, unlike $h$, is not well-defined modulo torsion. Recent attempts to find topologically better-behaved objects related to the Mahler measure include the introduction of the \emph{metric Mahler measure}, a well-defined metric on $\cF$, by Dubickas and Smyth \cite{DS}, and later the introduction of the \emph{ultrametric Mahler measure} by the first author and Samuels \cite{FS}. Both metrics induce the discrete topology if (and only if) Lehmer's conjecture is true.

In order to construct our norms related to the Mahler measure, we first construct an orthogonal decomposition of the space $\cF$ of algebraic numbers modulo torsion. We fix our algebraic closure $\Qbar$ of $\bQ$ and let $\cK$ denote the set of finite extensions of $\bQ$. We let $G=\Gal(\Qbar/\bQ)$ be the absolute Galois group, and let $\cK^G = \{ K\in\cK : \sigma K = K\text{ for all }\sigma\in G\}$. Let $V_K$ denote the $\bQ$-vector space span of the functions given by
\[
 V_K = \operatorname{span}_{\bQ}\langle \{\fal : \al\in K^\times/\Tor(K^\times) \}\rangle.
\]
We first prove the following result, which gives the orthogonal decomposition by Galois field:
\begin{thm}\label{thm:decomp-K}
 There exist projection operators $T_K :\cF\ra \cF$ for each $K\in\cK^G$ such that $T_K(\cF)\subset V_K$, $T_K(\cF)\perp T_L(\cF)$ for all $K\neq L\in\cK^G$ with respect to the inner product on $\cF$, and
 \[
  \cF = \bigoplus_{K\in\cK^G} T_K(\cF).
 \]
\end{thm}
In particular, we see that the projection operators $T_K$ are orthogonal projections with respect to the inner product on $\cF$, and thus in the completion with respect to the $L^2$ norm this gives a Hilbert space decomposition. A decomposition by Galois field alone, however, does not give enough information about the degree of a specific number in order to bound the Mahler measure of the number (and further, as we will see in Remark \ref{rmk:why-not-decompose-K}, a canonical decomposition along the entire collection of number fields is not possible). We therefore define the vector subspace
\[
 V^{(n)} = \sum_{\substack{K\in\cK\\ [K:\bQ]\leq n}} V_K
\]
and determine the following decomposition:
\begin{thm}\label{thm:decomp-n}
 There exist projections $T^{(n)}:\cF\ra\cF$ for each $n\in\bN$ such that $T^{(n)}(\cF)\subset V^{(n)}$, $T^{(m)}(\cF)\perp T^{(n)}(\cF)$ for all $m\neq n$, and
 \[
  \cF = \bigoplus_{n=1}^\infty T^{(n)}(\cF).
 \]
\end{thm}
These decompositions are independent of each other. Specifically, we have the following theorem:
\begin{thm}\label{thm:TK-Tn-commute}
 The projections $T_K$ and $T^{(n)}$ commute with each other for each $K\in\cK^G$ and $n\in\bN$.
\end{thm}
\noindent In other words, as a result of commutativity, we can form projections $T^{(n)}_K = T_K T^{(n)}$ and so we have an orthogonal decomposition
 \[
  \cF = \bigoplus_{n=1}^\infty \bigoplus_{K\in\cK^G} T^{(n)}_K(\cF).  
 \]
Again, when we pass to the completion in the $L^2$ norm, the projections extend by continuity and the above decomposition extends to the respective closures and the direct sum becomes a direct sum in the usual Hilbert space sense.

This geometric structure within the algebraic numbers allows us to define linear operators, for all $L^p$ norms with $1\leq p\leq \infty$, which capture the contribution of the degree to the Mahler measure in such a way that we can define our Mahler norms. Specifically, we define the operator
\begin{equation*}
\begin{split}
 M : \cF &\ra \cF\\
 f &\mapsto \sum_{n=1}^\infty n\,T^{(n)} f.
\end{split}
\end{equation*}
The sum is finite for each $f\in\cF$. $M$ is a well-defined, unbounded (in any $L^p$ norm, $1\leq p\leq \infty$), invertible linear map defined on the incomplete vector space $\cF$. We define the \emph{Mahler $p$-norm} on $\cF$ for $1\leq p\leq \infty$ to be
\begin{equation*} 
  \|f\|_{m,p} = \| Mf\|_p 
\end{equation*}
where $\| \cdot \|_p $ denotes the usual $L^p$ norm on the incomplete vector space $\cF$. The Mahler $p$-norm is, in fact, a well-defined vector space norm on $\cF$, and hence the completion $\cF_{m,p}$ with respect to $\|\cdot\|_{m,p}$ is a Banach space.

In order to see that these norms form a suitable generalization of the Mahler measure of algebraic numbers, we will show that the Lehmer conjecture can be equivalently reformulated in terms of these norms. First, let us address what form the Lehmer conjecture takes inside $\cF$. For any $\al\in\Qbar^\times$, let $h_p(\al)=\|\fal\|_p$. (Recall that $h_1(\al)=2\,h(\al)$.) Then we formulate:
\begin{conj}[$L^p$ Lehmer conjectures]\label{conj:lehmer-p}
 For $1\leq p\leq \infty$, there exists an absolute constant $c_p$ such that the $L^p$ Mahler measure satisfies the following equation:
 \begin{equation}\tag{$*_p$}\label{eqn:Lp-lehmer-conj}
 m_p(\al)=(\deg_{\bQ}\al)\cdot h_p(\al)\geq c_p>0\quad\text{for all}\quad \al\in\Qbar^\times\setminus \Tor(\Qbar^\times).
 \end{equation}
\end{conj}
\noindent From the fact that $h_1(\al)=2h(\al)$ it is clear that when $p=1$ this statement is equivalent to the Lehmer conjecture. For $p=\infty$, we will show in Proposition \ref{prop:equivalent-to-classical-conjectures} below that the statement is equivalent to the Schinzel-Zassenhaus conjecture.

In order to translate the Lehmer conjecture into a bound on function space norms which, unlike the metric Mahler measure, cannot possibly be discrete, it is necessary to reduce the Lehmer problem to a sufficiently small set of numbers which we can expect to be bounded away from zero in norm. This requires the introduction in Section \ref{sect:d-and-l} of two classes of algebraic numbers modulo torsion in $\cF$, the \emph{Lehmer irreducible} elements $\cL$ and the \emph{projection irreducible elements} $\cP$. Let $\cU\subset\cF$ denote the subspace of algebraic units. Then we prove the following theorem:
\begin{thm}\label{thm:Lp-Lehmer-equiv}
 For each $1\leq p\leq \infty$, equation \eqref{eqn:Lp-lehmer-conj} holds if and only if
 \begin{equation}\tag{$**_p$}\label{eqn:norm-Lp-lehmer}
  \|f\|_{m,p} \geq c_p > 0\quad\text{for all}\quad 0\neq f\in\cL\cap \cP\cap\cU
 \end{equation}
 where $\cL$ denotes the set of Lehmer irreducible elements, $\cP$ the set of projection irreducible elements, and $\cU$ the subspace of algebraic units. Further, for $1\leq p\leq q\leq\infty$, if ($**_p$) holds then ($**_q$) holds as well.
\end{thm}
\noindent The last statement of the theorem, which is proven by reducing to a place of measure 1 and applying the usual inequality for the $L^p$ and $L^q$ norms on a probability space, generalizes the well-known fact that Lehmer's conjecture implies the conjecture of Schinzel-Zassenhaus.

Let $\cU_{m,p}$ denote the Banach space which is the completion of the vector space $\cU$ of units with respect to the Mahler $p$-norm $\|\cdot\|_{m,p}$. The set $\cL\cap\cP\cap\cU$ has another useful property which we will prove, namely, that the additive subgroup it generates
\[
 \Gamma = \langle \cL\cap\cP\cap\cU\rangle
\]
is also a set of equivalence for the Lehmer conjecture, that is, we will show the $L^p$ Lehmer conjecture \eqref{eqn:Lp-lehmer-conj} is equivalent to the condition that $\Gamma$ be a discrete subgroup in $\cU_{m,p}$. Specifically, we have:
\begin{thm}\label{thm:Gamma-equiv}
 Equation \eqref{eqn:Lp-lehmer-conj} holds if and only if $\Gamma\subset \cU_{m,p}$ is closed.
\end{thm}
\noindent This result follows from the general fact that for a separable Banach space, an additive subgroup is discrete if and only it is closed and free abelian, and we will show that $\Gamma$ is in fact a free abelian group. This leads us to a new conjecture, equivalent to \eqref{eqn:Lp-lehmer-conj} for each $1\leq p\leq\infty$:
\begin{conj}\label{conj:closure-of-gamma}
 The group $\Gamma\subset\cU_{m,p}$ is closed for each $1\leq p\leq\infty$.
\end{conj}

Lastly, the presence of orthogonal decompositions raises a particular interest in the study of the $L^2$ norm. In this case, the norm associated to the Mahler measure has a particularly simple form which is in sympathy with the geometry of $L^2$.
\begin{thm}\label{thm:mahler-2-norm}
 The Mahler $2$-norm satisfies
 \[
  \|f\|_{m,2}^2 = \sum_{n=1}^\infty n^2\,\|T^{(n)}(f)\|_2^2= \sum_{K\in\cK^G}\sum_{n=1}^\infty n^2\,\|T^{(n)}_K(f)\|_2^2.
 \]
 Further, the Mahler $2$-norm arises from the inner product
 \[
  \langle f,g\rangle_m = \langle Mf,Mg\rangle = \sum_{n=1}^\infty n^2\, \langle T^{(n)} f,T^{(n)}g\rangle
  = \sum_{K\in\cK^G} \sum_{n=1}^\infty n^2\, \langle T_K^{(n)} f,T_K^{(n)}g\rangle
 \]
where $\langle f,g\rangle =\int_Y fg\,d\lambda$ denotes the usual inner product in $L^2(Y)$, and therefore the completion $\cF_{m,2}$ of $\cF$ with respect to the Mahler $2$-norm is a Hilbert space.
\end{thm}

The structure of this paper is as follows. In Section \ref{sect:decompositions} we introduce the basic operators and subspaces of our study, namely, those arising naturally from number fields and Galois isomorphisms. The proofs of Theorems \ref{thm:decomp-K}, \ref{thm:decomp-n} and \ref{thm:TK-Tn-commute} regarding the orthogonal decompositions of the space $\cF$ with respect to Galois field and degree will then be carried out in \ref{sect:sub-decomp-1}, \ref{sect:sub-decomp-2}, and \ref{sect:sub-decomp-3}. In Section \ref{sect:d-and-l} we prove our results regarding the reduction of the classical Lehmer problem and introduce the relevant classes of algebraic numbers which are essential to our theorems. Finally in Section \ref{sect:mahler-p} we introduce the Mahler $p$-norms and prove the remaining results.

\section{Orthogonal Decompositions}\label{sect:decompositions}

\subsection{Galois isometries}
Let $\cF_p$ denote the completion of $\cF$ with respect to the $L^p$ norm. By \cite[Theorems 1-3]{AV},
\[
 \cF_p = \begin{cases}
 \{f\in L^1(Y,\lambda) : \int_Y f\,d\lambda =0\} & \text{if }p=1\\
 L^p(Y,\lambda) & \text{if }1<p<\infty\\
 C_0(Y,\lambda) & \text{if }p=\infty.
 \end{cases}
\]
We begin by introducing our first class of operators, the isometries arising from Galois automorphisms. Let us recall how the Galois group acts on the places of an arbitrary Galois extension $K$. Suppose $\al\in K$, $v\in M_K$ is a place of $K$, and $\sigma\in G$. We define $\sigma v$ to be the place of $K$ given by $\|\al\|_{\sigma v} = \|\sigma^{-1}\al\|_{v}$, or in other words, $\|\sigma\al\|_{v}=\|\al\|_{\sigma^{-1}v}$.

\begin{lemma}\label{lemma:sigma-measure-preserving}
 Each $\sigma\in G$ is a measure-preserving homeomorphism of $(Y,\lambda)$.
\end{lemma}
\begin{proof}
 Recall from \cite{AV} that $Y = \varprojlim_K Y_K$ where $K$ ranges over the finite Galois extensions of $\bQ$ and $Y_K$ denotes the set of places of $K$ endowed with the discrete topology. That $\sigma : Y\ra Y$ is a well-defined bijection follows from the fact that $G$ gives a group action. Continuity of $\sigma$ and $\sigma^{-1}$ follow from \cite[Lemma 3]{AV}. It remains to show that $\sigma$ is measure-preserving, but this follows immediately from \cite[(4.6)]{AV}.
\end{proof}

In accordance with the action on places, we define for $\sigma\in G$ the operator
\[
  L_\sigma : \cF_p \ra \cF_p
\]
given by
\[
 (L_\sigma f)(y) = f(\sigma^{-1}y).
\]
Thus for $f_\al\in\cF$, we have $L_\sigma f_\al = f_{\sigma\al}$, and in particular $L_\sigma(\cF)\subseteq \cF$ for all $\sigma\in G$. Further, by our definition of the action on places, we have $L_\sigma L_\tau = L_{\sigma\tau}$.

Let $\cB(\cF_p)$ denote the bounded linear maps from $\cF_p$ to itself, and let $\cI(\cF_p)\subset \cB(\cF_p)$ denote the subgroup of isometries of $\cF_p$. 
By the construction of $\lambda$, each $\sigma\in G$ is a measure-preserving topological homeomorphism of the space of places $Y$, so it follows immediately that $L_\sigma$ is an isometry for all $1\leq p\leq \infty$, that is, $\|L_\sigma f\|_p = \|f\|_p$ for all $\sigma\in G$. Thus we have a natural map
\begin{align*}
 \rho : G & \ra \cI(\cF_p)\\
 \sigma &\mapsto L_\sigma
\end{align*}
where $(L_\sigma f)(y) = f(\sigma^{-1} y)$. We will show that $\rho$ gives an injective infinite-dimensional representation of the absolute Galois group (which is unitary in the case of $L^2$), and further, that the map $\rho$ is continuous if $G$ is endowed with its natural profinite topology and $\cI$ is endowed with the strong operator topology inherited from $\cB(\cF_p)$. (Recall that the strong operator topology, which is weaker than the norm topology, is defined as the weakest topology such that the evaluation maps $A \mapsto \|Af\|_p$ are continuous for every $f\in L^p$.)
\begin{prop}
 The map $\rho : G\ra \cI$ is injective, and it is continuous if $\cI$ is endowed with the strong operator topology and $G$ has the usual profinite topology.
\end{prop}
\begin{proof}
 First we will observe that the image $\rho(G)$ is discrete in the norm topology, so that $\rho$ is injective. To see this, fix $\sigma\neq \tau\in G$, so that there exists some finite Galois extension $K$ and an element $\al\in K^\times$ such that $\sigma\al\neq \tau\al$. By \cite[Theorem 3]{D}, we can find a rational integer $n$ such that $\beta = n+\al$ is torsion-free, that is, if $\beta/\beta'\neq 1$ then $\beta/\beta'\not\in\Tor(\Qbar^\times)$ for any conjugate $\beta'$ of $\beta$, and in particular, the conjugates of $\beta$ give rise to distinct functions in $\cF$. Thus $\sigma\beta\neq \tau\beta$ implies that $L_\sigma f_\beta \neq L_\tau f_\beta$, so in particular, there exists some place $v$ of $K$ such that $\sigma(Y(K,v))\neq \tau(Y(K,v))$ and are therefore disjoint sets. Choose a Galois extension $L/K$ with distinct places $w_1,w_2|v$. Since $L/K$ is Galois, the local degrees agree and so $\lambda(Y(L,w_1))=\lambda(Y(L,w_2))$ by \cite[Theorem 5]{AV}. Define
 \[
  f(y) = \begin{cases}
         1 & \text{if }y\in Y(L,w_1)\\
        -1 & \text{if }y\in Y(L,w_2)\\
	 0 & \text{otherwise.}
        \end{cases}
 \]
 Clearly $f\in \cF_p$ for all $1\leq p\leq \infty$ and $L_\sigma f$ and $L_\tau f$ have disjoint support. Thus,
 \[
  \|(L_\sigma -L_\tau)f\|_p = \|L_\sigma f\|_p +\|L_\tau f\|_p = 2\|f\|_p.
 \]
 But this implies that $2\leq \|L_\sigma - L_\tau\|\leq \|L_\sigma\| + \|L_\tau\|=2$ so $\|L_\sigma - L_\tau\|=2$. Thus the image $\rho(G)$ is discrete in the norm topology of $\cI$, and $\rho$ is injective.
 
 Let us now prove continuity. Recall that a basis for the strong operator topology on $\cI$ is given by sets of the form
 \[
  U = \{ A\in\cI : \|(A-B)f_i\|<\ep\text{ for all }1\leq i\leq k\}
 \]
 where $B\in\cI$, $f_1,\ldots,f_k$ is a finite set of functions in $L^p$, and $\ep>0$. Fix such an open set $U$ for a given $B=L_\sigma$ for some $\sigma\in G$. Approximate each $f_i$ by an element $g_i\in\cF$ such that $\|f_i-g_i\|_p<\ep/2$. Let $V_K$ be a subspace of $\cF$ containing $g_1,\ldots,g_k$. Let
 \[
  N = \{\tau \in G : \sigma|_K = \tau|_K\}.
 \]
 Then $N$ is an open subset of $G$ in the profinite topology. We claim that $\rho(N)\subseteq U$, and thus that $\rho$ is continuous. To see this, observe that for $\tau\in N$,
 \begin{multline*}
  \|(L_\tau-L_\sigma)f_i\|_p \leq 
  \|(L_\tau-L_\sigma)g_i\|_p + \|(L_\tau-L_\sigma)(f_i-g_i)\|_p\\
  < \|(L_\tau-L_\sigma)g_i\|_p + 2\cdot \ep/2=\ep
 \end{multline*}
 where $\|(L_\tau-L_\sigma)g_i\|_p=0$ because $g_i\in V_K$, and thus is locally constant on the sets $Y(K,v)$ for $v$ a place of $K$, and $\tau\in N$ implies that $\sigma$ and $\tau$ agree on $K$, so $L_\tau g_i = L_\sigma g_i$.
\end{proof}

\subsection{Subspaces associated to number fields}
We will now prove some lemmas regarding the relationship between the spaces $V_K$ and the Galois group. As in the introduction, let us define
\[
 \cK = \{ K/\bQ : [K:\bQ]<\infty\}\quad\text{and}\quad 
 \cK^G = \{ K\in\cK : \sigma K =K\ \forall \sigma\in G\}.
\]
As we shall have occasion to use them, let us recall the combinatorial properties of the sets $\cK$ and $\cK^G$ partially ordered by inclusion. Recall that $\cK$ and $\cK^G$ are \emph{lattices}, that is, partially ordered sets for which any two elements have a unique greatest lower bound, called the \emph{meet}, and a least upper bound, called the \emph{join}. Specfically, for any two fields $K,L$, the meet $K\wedge L$ is given by $K\cap L$ and the join $K\vee L$ is given by $KL$. If $K,L$ are Galois then both the meet (the intersection) and the join (the compositum) are Galois as well, thus $\cK^G$ is a lattice as well. Both lattices have a minimal element, namely $\bQ$, and are \emph{locally finite}, that is, between any two fixed elements we have a finite number of intermediate elements.

For each $K\in\cK$, let
\begin{equation}
 V_K = \operatorname{span}_{\bQ}\langle \{\fal : \al\in K^\times/\Tor(K^\times) \}\rangle.
\end{equation}
Then $V_K$ is the subspace of $\cF$ spanned by the functions arising from numbers of $K$. Suppose we fix an algebraic number $f\in\cF$. Then the set 
\[
\{ K\in\cK : f\in V_K\}
\]
forms a sublattice of $\cK$, and by the finiteness properties of $\cK$ this set must contain a unique minimal element.
\begin{defn}
 For any $f\in \cF$, the \emph{minimal field} is defined to be the minimal element of the set $\{ K\in\cK : f\in V_K\}$. We denote the minimal field of $f$ by $K_f$.
\end{defn}

\begin{lemma}
 For any $f\in\cF$, we have $\Stab_{G}(f) = \Gal(\Qbar/K_f) \leq G$.
\end{lemma}
\begin{notation}
 By $\Stab_G(f)$ we mean the $\sigma\in G$ such that $L_\sigma f = f$. As this tacit identification is convenient we shall use it throughout without further comment.
\end{notation}
\begin{proof}
 Let $f=\fal$. Then clearly $\Gal(\Qbar/K_f)\leq \Stab_G(f)$, as $\al^\ell\in K_f$ for some $\ell\in\bN$ by definition of $V_{K_f}$. To see the reverse implication, merely observe that $K_f=\bQ(\al^\ell)$ for some $\ell\in\bN$, as otherwise, there would be a proper subfield of $K_f$ which contains a power of $\al$, contradicting the definition of $K_f$.
\end{proof}
\begin{rmk}
 The minimal such exponent $\ell$ used above can in fact be uniquely associated to $f\in\cF$ and this will be vital to the concept of Lehmer irreducibility developed in Section \ref{sect:d-and-l} below.
\end{rmk}

\begin{lemma}\label{lemma:when-is-al-in-VK}
 For a given $f\in\cF$, we have $f\in V_K$ if and only if $L_\sigma f= f$ for all $\sigma\in \Gal(\Qbar/K)$.
\end{lemma}
\begin{proof}
 Necessity is obvious. To see that the condition is sufficient, observe that by definition of $K_f$, we have $f \in V_K$ if and only $K_f\subseteq K$, which is equivalent to $\Gal(\Qbar/K)\leq \Gal(\Qbar/K_f)$ under the Galois correspondence. But by the above lemma, $\Gal(\Qbar/K_f)=\Stab_G(f)$.
\end{proof}

\begin{prop}\label{prop:VK-is-unique-to-K}
 If $E,F\in\cK$, then we have $E\neq F$ if and only if $V_E\neq V_{F}$.
\end{prop}
\begin{proof}
  Suppose $E\neq F$ but $V_E=V_F$. Let $E=\bQ(\al)$. By \cite[Theorem 3]{D} we can find a rational integer $n$ such that $\beta=n+\al$ is torsion-free, that is, if $\beta/\beta'\neq 1$ then $\beta/\beta'\not\in\Tor(\Qbar^\times)$ for any conjugate $\beta'$ of $\beta$, and in particular, the conjugates of $\beta$ give rise to distinct functions in $\cF$. Observe therefore that $E=\bQ(\beta)$ and $\Stab_G(f_\beta)=\Gal(\Qbar/E)$. By the above if $f_\beta\in V_F$ then we must have $\Gal(\Qbar/F)\leq \Gal(\Qbar/E)$, or $E\subseteq F$. Repeating the same argument for a generator of $F$, we find that $F\subseteq E$ so $E=F$, a contradiction. The reverse implication is obvious.
\end{proof}
\begin{rmk}
 The above corollary is no longer true if we restrict our attention to the space of units $\cU\subset\cF$. This follows from the well known fact that CM extensions (totally imaginary quadratic extensions of totally real fields) have the same unit group modulo torsion as their base fields, the simplest example being $\bQ(i)/\bQ$.
\end{rmk}

\subsection{Orthogonal projections associated to number fields}
For $K\in\cK$, define the map $P_K : \cF \ra  V_K$ via
\[
  (P_Kf)(y) =\int_{H_K} (L_{\sigma}f)(y)d\nu(\sigma) 
\]
where $H_K=\Gal(\Qbar/K)$ and $\nu$ is the normalized (measure $1$) Haar measure of $H_K$. (Observe that, like $G$, $H_K$ is profinite and thus compact and possesses a Haar measure.) Let us prove that the map is well-defined. Since $f\in\cF$, it has a finite Galois orbit and thus a finite orbit under $H_K$. Let us partition $H_K$ into the $k = [H_K : \Stab_{H_K}(f)]$ cosets of equal measure by the translation invariance of the Haar measure. Denote these cosets by $\Stab_{H_K}(f)\sigma_1, \ldots, \Stab_{H_K}(f)\sigma_k$. Then
\[ 
P_K(f) = \frac{1}{k}\left( L_{\sigma_1} f + \cdots + L_{\sigma_k} f \right).
\]
But each $L_{\sigma_i} f\in\cF$ since $\cF$ is closed under the action of the Galois isometries. Thus if $f=\fal$, we have $L_{\sigma_i} f = f_{\sigma_i\al}$. Since $\cF$ is a vector space, $P_K(f)\in \cF$ as well. Further, it is stable under the action of $H_K$, and thus, by Lemma \ref{lemma:when-is-al-in-VK}, we have $P_K(f)\in V_K$. The map $P_K$ is in fact nothing more than the familiar algebraic norm down to $K$, subject to an appropriate normalization, that is, if $f_\beta = P_K \fal$, then we have 
\[
 \beta \equiv \left(\ovrln{\N^{K(\al)}_K \al}\right)^{1/[K(\al):K]}\mod\Tor(\Qbar^\times).
\]
(We note in passing that the norm map $N^{K(\al)}_K : K(\al)^\times \ra K^\times$ is well-defined modulo torsion.)

The following alternative formulation will also be helpful:
\begin{lemma}\label{lemma:PK-as-conditional-expectation}
 Let $K\in\cK$ and let $M_K$ denote the places of $K$. For each $v\in M_K$, let $\chi_v(y)$ be the characteristic function of the set $Y(K,v)$. Then
 \[
  P_K f(y) = \sum_{v\in M_K} \left(\frac{1}{\lambda(Y(K,v))}
  \int_{Y(K,v)} f(z)\,d\lambda(z)\right) \chi_v(y).
 \]
\end{lemma}
\noindent In other words, $P_K$ is essentially the conditional expectation with respect to the Borel $\sigma$-algebra generated by the set $\{Y(K,v) : v\in M_K\}$. Of course, $Y$ has infinite measure so this is not a conditional expectation in the usual sense from probability theory, although it shares many of the same properties. If we restrict to the space of units, that is, functions supported on the measure one space $Y(\bQ,\infty)$, then the restriction of $P_K$ to this space is indeed a conditional expectation.
\begin{proof}
 Fix a value $y\in Y$. Then there exists a unique $v\in M_K$ such that $y\in Y(K,v)$ since $Y=\bigcup_{v\in M_K} Y(K,v)$ is a disjoint union. The claim will be proven if we can show that for this value of $y$,
 \[
  P_K f(y) = \frac{1}{\lambda(Y(K,v))}
  \int_{Y(K,v)} f(z)\,d\lambda(z).
 \]
 Now,  
 \[
  P_K f(y) = \int_{H_K} f(\sigma^{-1}y)\,d\nu(\sigma)
 \]
 where $H_K,\nu$ are as above. By the construction of $\lambda$ (see (4.1) and surrounding remarks in \cite{AV}), for any $y\in Y(K,v)$,
 \[
  \frac{1}{\lambda(Y(K,v))}\int_{Y(K,v)} f(z)\,d\lambda(z) = \int_{H_K} f(\sigma^{-1}y)\,d\nu(\sigma)
 \]
 (where we need the normalization factor ${1}/{\lambda(Y(K,v))}$ since (4.1) assumes $\lambda(K,v)=1$) and so the proof is complete.
\end{proof}

\begin{prop}\label{prop:PK-continuous-wrt-Weil-norm}
 Let $\bQ\subset K\subset\Qbar$ be an arbitrary field. Then $P_K$ is a projection onto $V_K$ of norm one with respect to the $L^p$ norms for $1\leq p\leq \infty$.
\end{prop}
\begin{proof}
 We first prove that $P_K^2 = P_K$. Let $H = H_K$ as above and $\nu$ the normalized Haar measure on $H$. Suppose that $\tau\in H$. Observe that
\[
 P_K(f)(\tau^{-1} y) = \int_H f(\sigma^{-1} \tau^{-1} y)d\nu(\sigma) = \int_{\tau H} f(\sigma^{-1} y)d\nu(\sigma) = P_K(f)(y)
\]
 since $\tau H = H$ for $\tau\in H$.
 Thus,
\[
 ({P_K}^2f)(y) = \int_H P_K f(\sigma^{-1} y)d\nu(\sigma) =\int_H P_K f(y) d\nu(\sigma) = P_K f(y),
\]
 or more succinctly, ${P_K}^2=P_K$. Since linearity is clear we will now prove that the operator norm $\|P_K\|=1$ in the $L^p$ norm in order to conclude that $P_K$ is a projection. If $p=\infty$, this is immediate, so let us assume that $1\leq p<\infty$. Let $f\in L^p(Y)$. Then first observe that since $\nu(H)=1$, Jensen's inequality implies
\[
 \int_H |f(\sigma^{-1}y)|\,d\nu(\sigma) \leq \left(\int_H |f(\sigma^{-1}y)|^p\,d\nu(\sigma)\right)^{1/p}.
\]
 Now let us consider the $L^p$ norm of $P_K f$:
 \begin{multline*}
 \|P_Kf\|_p= \left(\int_Y |P_K(f)(y)|^p d\lambda(y)\right)^{1/p} = \left(\int_Y\left|\int_H f(\sigma^{-1} y)d\nu(\sigma)\right|^p d\lambda(y)\right)^{1/p}\\
 \leq \left( \int_Y \int_H \left|f(\sigma^{-1} y)\right|^p d\nu(\sigma)d\lambda(y) \right)^{1/p}=
 \left( \int_H \int_Y\left|f(\sigma^{-1} y)\right|^p d\lambda(y)d\mu(\sigma)\right)^{1/p}\\
 =\left(\int_H \|L_\sigma f\|_p^p d\mu(\sigma)\right)^{1/p}=\left(\int_H \|f\|_p^p d\mu(\sigma)\right)^{1/p}=\|f\|_p.
 \end{multline*}
 where we have made use of the fact that $L_\sigma$ is an isometry, and the application of Fubini's theorem is justified by the integrability of $|f|^p$. This proves that $\|P_K\|\leq 1$, and to see that the operator norm is not in fact less than $1$, observe that the subspace $V_{\bQ}$ is fixed for every $P_K$.
\end{proof}
As a corollary, if we extend $P_K$ by continuity to the completion $\cF_p$ of $\cF$ under the $L^p$ norm, we obtain:
\begin{cor}
 The subspace $\overline{V_K}\subset\cF_p$ is complemented in $\cF_p$ for all $1\leq p\leq \infty$.
\end{cor}

As $\cF_2=L^2(Y,\lambda)$ is a Hilbert space, more is in fact true:
\begin{prop}\label{prop:PK-is-orthogonal}
 For each $K\in\cK$, $P_K$ is the orthogonal projection onto the subspace $\overline{V_K}\subset L^2(Y)$.
\end{prop}
\noindent Specifically, this means that $\|f\|_2^2 = \|P_K f\|_2^2 + \|(I-P_K)f\|_2^2$, where $I$ is the identity operator.
\begin{proof}
 It suffices to observe that $P_K$ is idempotent and has operator norm $\|P_K\|=1$ with respect to the $L^2$ norm, and any such projection in a real Hilbert space is orthogonal.
\end{proof}

We now explore the relationship between the Galois isometries and the projection operators $P_K$ for $K\in\cK$.

\begin{lemma}\label{lemma:PKCommWithG}
 For any field $K\subseteq\Qbar$ and $\sigma\in G$, 
 \[
  L_\sigma P_K = P_{\sigma K}\, L_\sigma.
 \]
 Equivalently, $P_K\, L_\sigma = L_\sigma P_{\sigma^{-1}K}$.
\end{lemma}
\begin{proof}
 We prove the first form, the second obviously being equivalent. By definition of $P_K$, letting $H=\Gal(\Qbar/K)$ and $\nu$ be the normalized Haar measure on $H$ such that $\nu(H)=1$,
\begin{align*}
 (L_\sigma P_K f)(y) = (P_K f)(\sigma^{-1}y) &= \int_H f(\tau^{-1} \sigma^{-1} y)\,d\nu(\tau) \\
 &= \int_H f(\sigma^{-1} \sigma \tau^{-1} \sigma^{-1} y)\,d\nu(\tau) \\
 &= \int_H f(\sigma^{-1}(\sigma\tau\sigma^{-1})^{-1} y)\,d\nu(\tau) \\
 &= \int_H (L_\sigma f)((\sigma\tau\sigma^{-1})^{-1} y)\,d\nu(\tau) \\
 &= \int_{\sigma H\sigma^{-1}} (L_\sigma f)(\tau^{-1} y)\,d\nu(\tau) \\
 & = P_{\sigma K}(L_\sigma f)(y).\qedhere
\end{align*}
\end{proof}

We will be particularly interested in the case where the projections $P_K,P_L$ commute with each other (and thus $P_K P_L$ is a projection to the intersection of their ranges). To that end, let us determine the intersection of two distinguished subspaces:
\begin{lemma}\label{lemma:intersection-of-VK-VL}
 Let $K,L\subset\Qbar$ be extensions of $\bQ$ of arbitrary degree. Then the intersection $V_K\cap V_L = V_{K\cap L}$.
\end{lemma}
\begin{proof}
 Simply observe that $\fal\in V_K$ if and only if $\al^n\in K$ for some $n\in\bN$, likewise, suppose $\al^m\in L$. Then $\al^{nm}\in K\cap L$, so $\fal\in V_{K\cap L}$. The reverse inclusion is obvious.
\end{proof}
\begin{lemma}\label{lemma:PKcommute}
 Suppose $K\in\cK$ and $L\in\cK^G$. Then $P_K$ and $P_L$ commute, that is,
 \[
  P_K P_L = P_{K\cap L} = P_L P_K.
 \]
 In particular, the family of operators $\{P_K : K\in\cK^G\}$ is commuting.
\end{lemma}
\begin{proof}
 It suffices to prove $P_K(V_L)\subset V_L$, as this will imply that $P_K(V_L)\subset V_K\cap V_L=V_{K\cap L}$ by the above lemma, and thus that $P_K P_L$ is itself a projection onto $V_{K\cap L}$, and thus $P_K P_L=P_{K\cap L}$. Since $P_{K\cap L}$ is an orthogonal projection, it is equal to its adjoint, and we find that $P_{K\cap L} = P_L P_K$ as well. To prove that $P_K(V_L)\subset V_L$, observe that for $f\in V_L$,
 \[
  P_K(f) = \frac{1}{k}(L_{\sigma_1} f +\cdots L_{\sigma_k} f)
 \]
 where the $\sigma_i$ are right coset representatives of $\Gal(\Qbar/L)\cap \Gal(\Qbar/K)$ in $\Gal(\Qbar/K)$. However, $L_\sigma(V_L)=V_L$ for $\sigma\in G$ since $L$ is Galois, and thus, $P_K(f)\in V_L$ as well. But $P_K(f)\in V_K$ by construction and the proof is complete.
\end{proof}

\subsection{Main decomposition theorem}\label{sect:sub-decomp-1}
We will now begin the proof of Theorems \ref{thm:decomp-K} and \ref{thm:decomp-n}, which state that we can orthogonally decompose the space $\cF$ of algebraic numbers modulo torsion by their Galois field and by their degree. These results will be derived from the following general decomposition theorem, which we will apply to $\cF$ in the next two sections.
\begin{thm}\label{thm:decomposition}
Let $V$ be a vector space over $\bQ$ with an inner product and suppose we have a family of subspaces $V_i\subset V$ together with projections $P_i$ indexed by a partially ordered set $I$ such that:
\begin{enumerate}
 \item \label{list:1-locally-finite}The index set $I$ has a unique minimal element, denoted $0\in I$, and $I$ is locally finite, that is, any interval $[i,j]=\{k\in I : i\leq k\leq j\}$ is of finite cardinality.
 \item \label{list:2-meet-semilattice}Any pair of elements $i,j\in I$ has a unique greatest lower bound, called the meet of $i$ and $j$, and denoted $i\wedge j$. (Such a poset $I$ is called a meet-semilattice.)
 \item \label{list:3-containment}$V_i\subseteq V_j$ if $i\leq j\in I$.
 \item \label{list:P-ortho}The projection map $P_i : V\ra V_i$ is orthogonal for all $i\in I$.
 \item \label{list:P-commute}For $i,j\in I$, $P_i P_j = P_j P_i = P_{i\wedge j}$, where $i\wedge j$ is the meet of $i$ and $j$.
 \item \label{list:all-v}$V=\sum_{i\in I} V_i$.
\end{enumerate}
Then there exist mutually orthogonal projections $T_i\leq P_i$ (that is, satisfying $T_i(V)\subseteq V_i$) which form an orthogonal decomposition of $V$:
\[
 V = \bigoplus_{i\in I} T_i(V),\quad\text{and}\quad T_i(V)\perp T_j(V)\text{ for 
all }i\neq j\in I.
\]
\end{thm}
We call $T_i$ the \emph{essential projection} associated to the space $V_i$, as it gives the subspace of $V_i$ which is unique to $V_i$ and no other subspace $V_j$ in the given family.
\begin{rmk}
 Theorem \ref{thm:decomposition} can be stated and proven almost identically if $V$ is a real Hilbert space rather than an incomplete vector space over $\bQ$, the only changes being that condition \eqref{list:all-v} is replaced with the condition that the closure of $\sum_{i\in I} V_i$ is $V$, the direct sum is understood in the usual Hilbert space sense, and the expansion of each $f$ into $\sum_{i\in I} T_i f$ is to be understood as a series rather than a finite sum. The construction of the $T_i$ operators and the orthogonality are proven in exactly the same manner, and indeed, we will make use of the fact that if we complete $V$, the decomposition extends by continuity to the completion in the usual Hilbert space sense. The theorem as stated here and as applied to $\cF$ is in fact a strictly stronger result than the statement it implies for the decomposition of $L^2(Y)$ as not only must such projections and such a decomposition exist, but this decomposition must also respect the underlying $\bQ$ vector space of algebraic numbers $\cF$ and map algebraic numbers to algebraic numbers.
\end{rmk}
Let us begin by recalling the background necessary to define our $T_i$ projections. Since $I$ is locally finite, it is a basic theorem in combinatorics that there exists a M\"obius function $\mu : I\times I\ra \bZ$, defined inductively by the requirements that $\mu(i,i)=1$ for all $i\in I$, $\mu(i,j)=0$ for all $i\not\leq j\in I$, and
$
 \sum_{i\leq j\leq k} \mu(i,j) = 0 
$
for all $i,k\in I$ (the sums are finite by the assumption that $I$ is locally finite). Since our set $I$ has a minimal element $0$ and is locally finite, we can sum over $i\leq j$ as well. The most basic result concerning the M\"obius function is \emph{M\"obius inversion}, which (in one of the several possible formulations) tells us that given two functions $f,g$ on $I$,
\[
 f(j) = \sum_{i\leq j} g(i)\quad\text{if and only if}\quad 
 g(j) = \sum_{i\leq j} \mu(i,j)\,f(i).
\]
In order that our $T_i$ capture the unique contribution of each subfield $V_i$, we would like our $T_i$ projections to satisfy the condition that:
\[
 P_j = \sum_{i\leq j} T_i.
\]
M\"obius inversion leads us to define the $T_i$ operators via the equation:
\begin{equation}\label{eqn:T-defn}
 T_j = \sum_{i\leq j} \mu(i,j) P_i.
\end{equation}
Since each of the above sums is finite and $\mu$ takes values in $\bZ$, we see that $T_j : V\ra V_j$ is well-defined. We will prove that $T_j$ is desired the family of projections.

\begin{lemma}\label{lemma:PiTj}
 Let the projections $P_i$ for $i\in I$ satisfy the conditions of Theorem \ref{thm:decomposition} and let $T_i$ be defined as above. Then for all $i,j\in I$,
 $
  P_i T_j = T_j P_i,
 $
 and 
 \[
 P_j T_i = \begin{cases}
               T_j & \text{if } i\leq j \\
               0 & \text{otherwise.} 
            \end{cases}
 \]
\end{lemma}
\begin{proof}
 The first claim follows immediately from equation \eqref{eqn:T-defn} and condition \eqref{list:P-commute} of the theorem statement. To prove the second claim, we proceed by induction.  Observe that the statement is trivial for $T_0 = P_0$. Now given $j\in I$, suppose the theorem is true for all $i<j$. Observe that from \eqref{eqn:T-defn} we get
 \begin{equation}\label{eqn:T-useful}
  T_j = P_j-\sum_{ i<j } T_i.  
 \end{equation}
 Then, if $i<j$, we have
 \[
  P_jT_i=P_jP_i-\sum_{k<i}P_jT_k=P_i-\sum_{k<i}T_k=T_i,
 \]
 applying the induction hypothesis at the second equality.

Now suppose $i\not<j$,  so that $i\wedge j\neq i$.  Then
 \begin{multline*}
  P_j T_i = P_jP_i-\sum_{k< i} P_jT_k = P_{i\wedge j}-\sum_{k\leq i\wedge j} P_jT_k - \sum_{\substack{k< i\\ k\not \leq i\wedge j}}P_jT_k\\
  = P_{i\wedge j}-\sum_{k\leq i\wedge j} T_k-0= P_{i\wedge j}-P_{i\wedge j}=0
 \end{multline*}
by two applications of the induction hypothesis at the third equality.
\end{proof}

\begin{lemma}\label{lemma:distinct-T}
 Let the $T_i$ be as above and let $i\neq j$ for $i,j\in I$. Then $T_i T_j = T_j T_i=0$.
\end{lemma}
\begin{proof}
 Suppose that $i\wedge j< j$. By Lemma \ref{lemma:PiTj}, $T_i = T_i P_i$ and $T_j = P_j T_j$. Thus,
 \[
  T_i T_j = (T_i P_i)(P_j T_j) = T_i(P_i P_j)T_j = T_i P_{i\wedge j} T_j = 0
 \]
 since $i\neq j$ implies that $i\wedge j< i$ or $i\wedge j< j$, so either $T_i P_{i\wedge j}=0$ or $P_{i\wedge j} T_j=0$ by Lemma \ref{lemma:PiTj}.
\end{proof}
We are now ready to prove the theorem statement.
\begin{proof}[Proof of Theorem \ref{thm:decomposition}]
 Let the operators $T_i$ for $i\in I$ be constructed as above. Let us first show that each $T_i$ is a projection, a linear operator of bounded norm such that ${T_i}^2 = T_i$. The fact the $T_i$ is a continuous linear operator of bounded norm follows from the same fact for the $P_i$ operators, since each $T_i$ is a finite linear combination of $P_i$ projections.
 
 Let us now show that $T_i$ is idempotent. The base case $T_0=P_0$ is trivial. Assume the lemma is true for all subfields for all $i<j$. Using equation \eqref{eqn:T-useful}, we have
 \begin{multline*}
  {T_j}^2 = \bigg(P_j-\sum_{ i< j} T_i\bigg)^2 = {P_j}^2 - \sum_{ i< j} P_j T_i - \sum_{ i< j} T_i P_j + \bigg(\sum_{ i< j} T_i\bigg)^2\\
  = {P_j} - \sum_{ i< j} T_i - \sum_{ i< j} T_i + \sum_{ i< j} T_i
  = {P_j} - \sum_{ i< j} T_i = T_j
 \end{multline*}
 where we have used Lemmas \ref{lemma:PiTj} and \ref{lemma:distinct-T} to simplify the middle and last terms.
 
 Now, let us show that the $T_i$ decompose $V$. To see this, observe that each element $f\in V$ by condition \eqref{list:all-v} lies in some $V_{i_1}+\ldots+V_{i_n}$. Let $I'=\bigcup_{m=1}^n [0,i_m]\subset I$, and then observe that $\sum_{k\in A} T_k$ is the projection onto $V_{i_1}+\ldots+V_{i_n}$ and $I'$ is finite by construction, so
 $
  f = \sum_{k\in I'} T_k f.
 $
 In fact, observe that we can write
 $
  f = \sum_{k\in I} T_k f
 $
 as a formally infinite sum, and all terms except those satisfying $k\leq i$ are zero by Lemma \ref{lemma:PiTj}. Thus we can write
 \[
  V = \bigoplus_{i\in I} T_i(V)
 \]
 and the fact that the $T_i$ are orthogonal projections now follows from this decomposition and Lemma \ref{lemma:distinct-T}.
\end{proof}

\subsection{Decomposition by Galois field and proof of Theorem \ref{thm:decomp-K}}\label{sect:sub-decomp-2}
We will now apply Theorem \ref{thm:decomposition} to $\cF$. Recall that $\cK^G$ is simply the set of finite Galois extensions of $\bQ$. As remarked above, it is well known that both $\cK$ and $\cK^G$ satisfy all of the axioms of a lattice, that is, for any two field $K,L$, there is a unique meet $K\wedge L$ given by $K\cap L$ and a unique \emph{join} $K\vee L$ given by $KL$. If $K,L$ are Galois then both the meet (the intersection) and the join (the compositum) are Galois as well, thus $\cK^G$ is a lattice as well. Further, both $\cK$ and $\cK^G$ are locally finite posets and possess a minimal element, namely, $\bQ$.

Our decomposition will be along $\cK^G$ and the associated family of subspaces $V_K$ with their canonical projections $P_K$. By since $\cK^G$ is a locally finite lattice, conditions \eqref{list:1-locally-finite} and \eqref{list:2-meet-semilattice} of Theorem \ref{thm:decomposition} are satisfied. Clearly the subspaces $V_K$ for $K\in\cK^G$ satisfy the containment condition \eqref{list:3-containment}. By Proposition \ref{prop:PK-is-orthogonal}, the projections are orthogonal and satisfy condition \eqref{list:P-ortho}. By Lemma \ref{lemma:PKcommute}, the maps $\{P_K : K\in\cK^G\}$ form a commuting family and satisfy condition \eqref{list:P-commute}. Lastly, since any $f=\fal$ belongs to $V_{K_f}\subset V_K$ where $K\in\cK^G$ is the Galois closure of the minimal field $K_f$, we find that condition \eqref{list:all-v} is satisfied as well. Thus Theorem \ref{thm:decomposition} gives us an orthgonal decomposition
\begin{equation}\label{eqn:decomposition-by-field}
 \cF = \bigoplus_{K\in\cK^G} T_K(\cF)
\end{equation}
The relationship between the $P_K$ and $T_K$ operators is given by:
\begin{equation}\label{eqn:TKdefn}
 P_K = \sum_{\substack{F\in\cK^G\\ F\subseteq K}} T_F,
 \quad\text{and}\quad T_K =\sum_{\substack{F\in\cK^G\\ F\subseteq K}} \mu(F,K)P_F 
\end{equation}
where $\mu : \cK^G\times \cK^G\ra\bZ$ is the M\"obius function associated to $\cK^G$.

If $K$ is the Galois closure of the minimal field $K_f$ where $f=\fal$, then $P_K(f) = f$, and so \eqref{eqn:TKdefn} gives us a unique representation modulo torsion of the algebraic number $\bal$ which we call the \emph{$M$-factorization} of $\bal$, or the \emph{$M$-expansion} of $\fal$ in functional notation.
\begin{example}\label{ex:t-factorization}
 Let $\bal = 2+\sqrt{2}$ and let $f=\fal$. Then $K_f = \bQ(\sqrt{2})$. Since $K\in\cK^G$, $[K:\bQ]=2$ and it is easy to see that the interval $[\bQ,K]=\{\bQ,K\}\subset \cK^G$, and so $\mu(\bQ,K)=-1$, and thus
 \[
  T_K = P_K - P_{\bQ},\quad T_{\bQ} = P_{\bQ}.
 \]
 Thus
 \[
  T_K(\fal) = f_{1+\sqrt{2}},\quad T_{\bQ}(\fal) = f_{\sqrt{2}},
 \]
 and the $M$-factorization of $\bal$ has the form
 $
  {2+\sqrt{2}} = {\sqrt{2}}\cdot {(1+\sqrt{2})},
 $
 or in functional notation,
 \[
  f_{2+\sqrt{2}} = f_{\sqrt{2}} + f_{1+\sqrt{2}},
  \quad\text{and}\quad f_{\sqrt{2}} \perp f_{1+\sqrt{2}}.
 \]
\end{example}

\begin{rmk}\label{rmk:why-not-decompose-K}
 We end this section with a remark on why we decompose along $\cK^G$ but not $\cK$. It is not difficult to see that the $P_K$ projections for $K\in\cK$ do not form a commuting family. To see this, suppose $\al$ is a cubic algebraic unit with conjugates $\beta,\gamma$ and discriminant $\Delta$. Then we have the following fields:
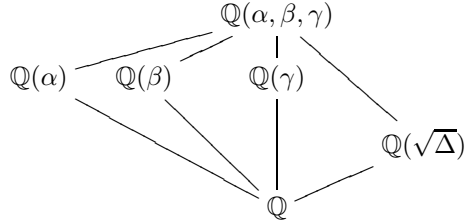
\[\xymatrix@R-=8pt@C=12pt{
 & & \bQ(\al,\beta,\gamma) & & \\
 \bQ(\al)\ar@{-}[urr] & \bQ(\beta)\ar@{-}[ur] & \bQ(\gamma)\ar@{-}[u] & &  \\
 & & & \bQ(\sqrt{\Delta})\ar@{-}[uul]\\
 & & \bQ\ar@{-}[uu]\ar@{-}[uul]\ar@{-}[uull]\ar@{-}[ur] & & 
 }
\]
 But the projections associated to the fields $\bQ(\al)$ and its conjugates do not commute. Specifically, we may compute:
\[
 P_{\bQ(\beta)}f_\al = -\frac{1}{2} f_\beta,\quad\text{and}\quad
 P_{\bQ(\al)}f_\beta = -\frac{1}{2} f_\al
\]
 which shows that $P_{\bQ(\al)}P_{\bQ(\beta)} \neq P_{\bQ(\beta)}P_{\bQ(\al)}$.  This noncommutativity is present precisely because there is a linear dependence among the vector space $V_{\bQ(\al)}$ and its conjugates, e.g., $f_\al+f_\beta+f_\gamma=0$ (since we assumed $\al$ was an algebraic unit). In particular, it is not hard to check that
 \[
  V_{\bQ(\al)}+V_{\bQ(\beta)}=V_{\bQ(\al)}+V_{\bQ(\beta)}+V_{\bQ(\gamma)}.
 \]
 Clearly such a dependence would make it impossible to associate a unique component $T_K$ to each of the three fields. However, the commutavity of the $P_K$ for $K\in\cK^G$ implies that there is no such barrier to decomposition amongst the Galois fields.
\end{rmk}

\subsection{Decomposition by degree and proof of Theorems \ref{thm:decomp-n} and \ref{thm:TK-Tn-commute}}\label{sect:sub-decomp-3}
In order to associate a notion of degree to a subspace in a meaningful fashion so that we can define our Mahler $p$-norms we will determine a decomposition of $\cF$. Let us define the function $\delta : \cF\ra \bN$ by
\begin{equation}\label{eqn:delta-defn}
 \delta(f) = \#\{L_\sigma f : \sigma\in G\}=[G:\Stab_G(f)]=[K_f:\bQ]
\end{equation}
to be the size of the orbit of $f$ under the action of the Galois isometries. Let 
\begin{equation}
 V^{(n)} = \sum_{\substack{K\in\cK\\ [K:\bQ]\leq n}} V_K
\end{equation}
be the vector space spanned by all elements of whose orbit in $\cF$ under $G$ is of size at most $n$. Let $P^{(n)}$ denote the orthogonal projection in $L^2(Y)$ onto the closure of $V^{(n)}$ in $L^2(Y)$. We wish to show that the restriction $P^{(n)}:\cF\ra V^{(n)}$ is a well-defined map of the algebraic numbers modulo torsion so that we can apply Theorem \ref{thm:decomposition} to construct projections $T^{(n)} : \cF\ra V^{(n)}$ which will give us the orthogonal decomposition of $\cF$ into a subspace spanned by elements whose orbit under $G$ is of size at most $n$. In order to prove this, we will first show that the projections $P^{(n)}$ and $P_K$ for $n\in\bN$ and $K\in\cK^G$ commute.

\begin{lemma}\label{lemma:delta-PK-leq-for-K-in-Kg}
 If $K\in\cK^G$, then $\delta(P_Kf)\leq \delta(f)$ for all $f\in\cF$.
\end{lemma}
\begin{proof}
 Let $F= K_f$. Since $K\in\cK^G$, we have by Lemma \ref{lemma:PKcommute} that $P_K f = P_K(P_F f) = P_{K\cap F} f$. Thus, $P_K f \in V_{K\cap F}$, and so $\delta(P_Kf)\leq [K\cap F:\bQ]\leq [F:\bQ]=\delta(f)$.
\end{proof}

\begin{prop}\label{prop:PnPK-commute}
 Let $n\in\bN$ and $K\in\cK^G$. Then the orthogonal projections $P^{(n)} : L^2(Y)\ra \overline{V^{(n)}}$ and $P_K : L^2(Y)\ra \overline{V_K}$ commute (where the closures are taken in $L^2$), and thus $T_K$ and $P^{(n)}$ commute as well.
\end{prop}
\begin{proof}
 Since $\delta(P_K f)\leq \delta(f)$ for all $f\in\cF$ by Lemma \ref{lemma:delta-PK-leq-for-K-in-Kg} above, we have $P_K(V^{(n)})\subset V^{(n)}$, and thus by continuity $P_K(\overline{V^{(n)}})\subset \overline{V^{(n)}}$, so $P_K(\overline{V^{(n)}})\subset \overline{V^{(n)}}\cap \overline{V_K}$ and $P_K P^{(n)}$ is a projection. Therefore they commute. The last part of the claim now follows from the definition of $T_K$ in \eqref{eqn:T-defn}.
\end{proof}

Let $W_K = T_K(\cF)\subset V_K$ for $K\in\cK^G$. By the above proposition, we see that if we can show that $P^{(n)}(W_K)\subseteq W_K$, then we will have the desired result, since 
\[
 P^{(n)}(\cF) = \bigoplus_{K\in\cK^G} P^{(n)}(W_K)
\]
by the commutativity of $P^{(n)}$ and $T_K$. Since we will prove this by reducing to finite dimensional $S$-unit subspaces, let us first prove an easy lemma regarding finite dimensional vector spaces over $\bQ$.
\begin{lemma}\label{lemma:sums-of-finite-diml-v-spaces}
 Suppose we have a finite dimensional vector space $A$ over $\bQ$, and suppose that
 \[
  A = V_1\oplus V_1' = V_2 \oplus V_2' = \cdots = V_n\oplus V_n'
 \]
 for some subspaces $V_i,V_i'$, $1\leq i\leq n$. Then
 \[
  A = (V_1 + \cdots + V_n)\oplus (V_1'\cap \cdots \cap V_n').
 \]
\end{lemma}
\begin{proof}
 It suffices to prove the lemma in the case $n=2$ as the remaining cases follow by induction, so suppose $A = V_1\oplus V_1' = V_2 \oplus V_2'$. It is an easy exercise that
 \[
  \dim_{\bQ} V_1 + \dim_{\bQ} V_2 = \dim_{\bQ} (V_1 + V_2) + \dim_{\bQ} (V_1\cap V_2),
 \]
 and likewise,
 \[
  \dim_{\bQ} V_1' + \dim_{\bQ} V_2' = \dim_{\bQ} (V_1' + V_2') + \dim_{\bQ} (V_1'\cap V_2').
 \]
 Now,
 \begin{multline}\label{eqn:temp}
  2\dim_{\bQ} A = \dim_{\bQ} V_1 + \dim_{\bQ} V_1' + \dim_{\bQ} V_2 + \dim_{\bQ} V_2'\\
  = \dim_{\bQ} (V_1 + V_2) + \dim_{\bQ} (V_1\cap V_2) + \dim_{\bQ} (V_1' + V_2') + \dim_{\bQ} (V_1'\cap V_2').
 \end{multline}
 Now, $(V_1 + V_2)\oplus (V_1'\cap V_2')\subseteq A$ and $(V_1' + V_2')\oplus (V_1\cap V_2)\subseteq A$, so
 \begin{align*}
  b &= \dim_{\bQ} (V_1 + V_2)+\dim_{\bQ} (V_1'\cap V_2')\leq \dim_{\bQ} A\\
  c &= \dim_{\bQ} (V_1' + V_2')+\dim_{\bQ} (V_1\cap V_2)\leq \dim_{\bQ} A.
 \end{align*}
 By \eqref{eqn:temp}, we have $b+c=2\dim_{\bQ} A$, therefore, we must have $b=c=\dim_{\bQ} A$, and in particular $b=\dim_{\bQ} A$ proves the claim.
\end{proof}

\begin{prop}
 $P^{(n)}(W_K)\subseteq W_K$ for every $n\in\bN$ and $K\in\cK^G$, and thus $P^{(n)}(\cF)\subset\cF$.
\end{prop}
\begin{proof}
 Let $f\in W_K$, and let $S\subset M_{\bQ}$ be a finite set of rational primes, containing the infinite prime, such that
 \[
  \supp_Y(f) \subset \bigcup_{p\in S} Y(\bQ,p).
 \]
 Let $V_{K,S}\subset V_K$ denote the subspace spanned by the $S$-units of $K$. By Dirichlet's $S$-unit theorem, $V_{K,S}$ is finite dimensional over $\bQ$. Let $W_{K,S}= T_K(V_{K,S})$. Notice that $W_{K,S}\subset V_{K,S}$ since each $P_F$ for each $F\subseteq K,\ F\in\cK^G$ will preserve the support of $f$ over each set $Y(\bQ,p)$ for $p\in M_{\bQ}$ by Lemma \ref{lemma:PK-as-conditional-expectation}. 
 
 For a field $F\in\cK$ such that $F\subset K$, let $W_{F,S} = P_F(W_{K,S})$ and $W_{F,S}' = Q_F(W_{K,S})$, where $Q_F = I - P_F$ is the complementary orthogonal projection. Observe that $W_{K,S} = W_{F,S} \oplus W_{F,S}'$. Then by Lemma \ref{lemma:sums-of-finite-diml-v-spaces}, we have
 \[
  W_{K,S} = \bigg(\sum_{\substack{F\subseteq K\\ [F:\bQ]\leq n}} W_{F,S}\bigg) \oplus \bigg(\bigcap_{\substack{F\subseteq K\\ [F:\bQ]\leq n}} W_{F,S}'\bigg).
 \]
 This gives us a decomposition $f=f_n + f_n'$ where
 \[
  f_n \in \sum_{\substack{F\subseteq K\\ [F:\bQ]\leq n}} W_{F,S} = V^{(n)} \cap V_{K,S},
 \]
 and
 \[
  f_n' \in \bigcap_{\substack{F\subseteq K\\ [F:\bQ]\leq n}} W_{F,S}'=(V^{(n)})^\perp \cap V_{K,S},
 \]
 But then $f_n\in V^{(n)}$ and $f_n'\in (V^{(n)})^\perp$, so by the uniqueness of the orthogonal decomposition, we must in fact have $f_n = P^{(n)} f$ and $f_n' = Q^{(n)}f = (I-P^{(n)})f$. Since this proof works for any $f\in\cF$, we have established the desired claim.
\end{proof}

Now we observe that the subspaces $V^{(n)}$ with their associated projections $P^{(n)}$, indexed by $\bN$ with the usual partial order $\leq$, satisfy the conditions of Theorem \ref{thm:decomposition}, and thus we have orthogonal projections $T^{(n)}$ and an orthogonal decomposition
\begin{equation}\label{eqn:decomposition-by-degree}
 \cF = \bigoplus_{n=1}^\infty T^{(n)}(\cF).
\end{equation}

The operators $T^{(n)}$ have a particularly simple form in terms of the $P^{(n)}$ projections. The M\"obius function for $\bN$ under the partial order $\leq$ is well-known and is merely 
\[
 \mu_{\bN}(m,n)=\begin{cases}
                 1 & \text{if }m=n,\\
		 -1 & \text{if }m=n-1\text{, and}\\
		 0 & \text{otherwise.}
                \end{cases}
\]
Thus, $T^{(1)}=P^{(1)}=P_{\bQ}$ and 
\[
 T^{(n)} = P^{(n)}-P^{(n-1)}\quad\text{for all}\quad n>1.
\]
We call $T^{(n)}f$ the \emph{degree $n$ component of $f$}. The following proposition is now obvious from the above constructions:
\begin{prop}\label{prop:Tn-expansion}
Each $f\in\cF$ has a unique finite expansion into its degree $n$ components, $f^{(n)}=T^{(n)}f\in\cF$
\[
 f = \sum_{n\in\bN} f^{(n)}.
\]
Each $f^{(n)}$ term can be written as a finite sum $f^{(n)}=\sum_i f^{(n)}_i$ where $f^{(n)}_i\in\cF$ and $\delta(f^{(n)}_i)=n$ for each $i$, and $f^{(n)}$ cannot be expressed as a finite sum $\sum_j f^{(n)}_j$ with $\delta(f^{(n)}_j)\leq n$ for each $j$ and $\delta(f^{(n)}_j)<n$ for some $j$.
\end{prop}
This completes the proof of Theorem \ref{thm:decomp-n}. It remains to prove Theorem \ref{thm:TK-Tn-commute}.
\begin{proof}[Proof of Theorem \ref{thm:TK-Tn-commute}]
 From Proposition \ref{prop:PnPK-commute}, we see that the operators $T_K$ and $P^{(n)}$ commute for $K\in\cK^G$ and $n\in\bN$. But $T^{(n)}=P^{(n)}-P^{(n-1)}$ for $n>1$ and $T^{(1)}=P^{(1)}$, so by the commutativity of $T_K$ with $P^{(n)}$ we have the desired result. In particular, the map $T^{(n)}_K = T^{(n)} T_K : \cF\ra\cF$ is also a projection, and thus we can combine equations \eqref{eqn:decomposition-by-field} and \eqref{eqn:decomposition-by-degree} to obtain the orthogonal decomposition
\begin{equation}\label{eqn:decomposition-by-deg-and-field}
 \cF = \bigoplus_{n=1}^\infty \bigoplus_{K\in\cK^G} T^{(n)}_K(\cF).\qedhere
\end{equation}
\end{proof}

\section{Reducing the Lehmer problem}\label{sect:d-and-l}
\subsection{Lehmer irreducibility}
Let us recall that we defined in Section \ref{sect:sub-decomp-3} the function $\delta : \cF\ra \bN$ by
\[
 \delta(f) = \#\{L_\sigma f : \sigma\in G\}=[G:\Stab_G(f)]=[K_f:\bQ].
\]
Observe that since nonzero scaling of $f$ does not affect its $\bQ$-vector space span or the minimal field $K_f$ that the function $\delta$ is invariant under nonzero scaling in $\cF$, that is, 
\[
\delta(rf)=\delta(f)\quad\text{for all}\quad  f\in\cF\text{ and }0\neq r\in\bQ.
\]
In order to better understand the relationship between our functions in $\cF$ and the algebraic numbers from which they arise, we need to understand when a function $\fal\in V_K$ has a representative $\al\in K^\times$ or merely is a root of an element $\al^n\in K^\times$ for some $n>1$. Naturally, the choice of coset representative modulo torsion affects this, and we would like to avoid such considerations. Therefore we define the function $d:\cF\ra\bN$ by
\begin{equation}\label{eqn:d-defn}
 d(f) = \min\{\deg_{\bQ} \al : \al\in\Qbar^\times,\ \fal=f\}.
\end{equation}
In other words, for a given function $f\in\cF$, which is an equivalence class of an algebraic number modulo torsion, $d(f)$ gives us the minimum degree amongst all of the coset representatives of $f$ in $\Qbar^\times$ modulo the torsion subgroup.

A number $f\in\cF$ can then be written as $f=\fal$ with $\al\in K_f^\times$ if and only if $d(f)=\delta(f)$. We therefore make the following definition:
\begin{defn}
 We define the set of \emph{Lehmer irreducible} elements of $\cF$ to be the set
\begin{equation}\label{eqn:LI-defn}
 \cL = \{f\in\cF : \delta(f) = d(f) \}.
\end{equation}
The set $\cL$ consists precisely of the functions $f$ such that $f=\fal$ for some $\al$ of degree equal to the degree of the minimal field of definition $K_f$ of $f$.
\end{defn}
We recall the terminology from \cite{D} that a number $\al\in\Qbar^\times$ is \emph{torsion-free} if $\al/\sigma\al\not\in\Tor(\Qbar^\times)$ for all distinct Galois conjugates $\sigma\al$. As we observed above in the proof of Proposition \ref{prop:VK-is-unique-to-K}, torsion-free numbers give rise to distinct functions $f_{\sigma\al}=L_\sigma \fal$ for each distinct Galois conjugate $\sigma\al$ of $\al$.
The goal of this subsection is to prove the following result relating $\delta$ and $d$:
\begin{prop}\label{prop:d-and-delta}
 Let $f\in\cF$ and $r,s\in\bZ$ with $(r,s)=1$. Then $R(f)=\{r\in\bQ : rf\in\cL\}=\frac{\ell}{n}\bZ$ where $\ell,n\in\bN$, $(\ell,n)=1$, and
 \begin{equation}
  d((r/s)f) = \frac{\ell s}{(\ell,r)(n,s)} \delta(f).
 \end{equation}
 In particular, $d(f)=\ell(f)\delta(f)$.
\end{prop}
We begin with the following lemma.
\begin{lemma}
 We have the following results:
 \begin{enumerate}
  \item For each $f\in\cF$, there is a unique minimal exponent $\ell(f)\in\bN$ such that $\ell(f) f\in\cL$.
  \item For any $\al\in\Qbar^\times$, we have $\delta(\fal) | \deg_{\bQ}\al$.
  \item $f\in\cL$ if and only if it has a representative in $\Qbar^\times$ which is torsion-free.
 \end{enumerate}
\end{lemma}
\begin{proof}
 Choose a representative $\al\in\Qbar^\times$ such that $f=\fal$ and let 
 \[
  \ell = \lcm\{ \ord(\al/\sigma\al) : \sigma\in G\text{ and }\al/\sigma\al\in\Tor(\Qbar^\times)\}
 \]
 where $\ord(\zeta)$ denotes the order of an element $\zeta\in\Tor(\Qbar^\times)$. Then observe that $\al^\ell$ is torsion-free. Clearly, $\bQ(\al^\ell)\subset \bQ(\al)$ so $[\bQ(\al^\ell):\bQ] | [\bQ(\al):\bQ]$. Now if a number $\beta\in\Qbar^\times$ is torsion-free, then since each distinct conjugate $\sigma\beta$ gives rise to a distinct function in $\cF$, we have
 \[
  \deg_{\bQ} \beta = [G : \Stab_G(f_\beta)] = [K_{f_\beta} : \bQ]=\delta(f_\beta).
 \]
 Thus $\deg_{\bQ} \al^\ell = \delta(\fal)$ and we have proven existence in the first claim. The existence of a minimum value follows since $\bN$ is discrete. To prove the second it now suffices to observe that since $\delta$ is invariant under scaling, with the choice of $\ell$ as above, we have $\delta(\fal) = \delta(\fal^\ell) | \deg_{\bQ} \al$ for all $\al\in\Qbar^\times$, and we have proven the second claim. The third now follows immediately.
\end{proof}
We note the following easy corollary for its independent interest:
\begin{cor}
 Let $\alpha\in\Qbar^\times$ have minimal polynomial $F(x)\in\bZ[x]$. Let $G(x)\in\bZ[x]$ be an irreducible polynomial of smallest degree in $\bZ[x]$ such that there exists some $k\in\bN$ with $F(x)|G(x^k)$. Then $\delta(\fal) = \deg G$.
\end{cor}
We note that $\delta(f)=1$ if and only if $f\in V_{\bQ}$, in which case, $f=\fal$ where $\al^n\in\bQ^\times$ and so $f$ represents a \emph{surd}.
\begin{lemma}\label{lemma:exponents-form-frac-ideal}
 If $0\neq f\in\cF$, then $R(f)=\{r\in\bQ : rf\in\cL\}$ is a fractional ideal of  $\bQ$.
\end{lemma}
\begin{proof}
 We can assume $\delta(f)>1$, otherwise $f$ arises from a surd and the proof is trivial. First we show that $R(f)$ is a $\bZ$-module. It is trivial that if $r\in R(f)$ then $-r\in R(f)$ as inversion does not affect degree. Suppose now that we have $r,s\in R(f)$ and choose torsion-free representatives $\beta^r,\gamma^s\in\Qbar^\times$ such that $\deg_{\bQ} \beta^r=\deg_{\bQ} \gamma^s=\delta(\bal)$. Since $\beta,\gamma$ both represent $f$, we have $\beta\equiv \gamma\mod\Tor(\Qbar^\times)$, we have $\beta=\zeta\gamma$ for some $\zeta\in\Tor(\Qbar^\times)$. Suppose $\zeta$ has order $N$. Then $\beta^{rN}$ and $\gamma^{sN}$ are both torsion-free and lie in the same field, and hence
 \[
 \bQ(\beta^r)=\bQ(\beta^{rN})=\bQ(\gamma^{sN})=\bQ(\gamma^s).
 \]
 But then $\beta^r\gamma^s= \zeta^r \gamma^{r+s}$ lies in this field as well. Since it is a representative of $(r+s)f$, it has degree at least $\delta(f)$ if $r+s\neq 0$ (which we can assume as otherwise the statement would be trivial) and $\gamma^{r+s}$ is nontorsion (otherwise $f$ itself would represent a torsion element and we would have $f=0$). Therefore, since it lies in a field of degree $\delta(f)$, it has degree $\delta(f)$. We conclude that $(r+s)f\in\cL$ as it has a representative of the requisite degree.
 
 If we can now show that $R(f)$ is finitely generated the proof will be complete.  But were it to require an infinite number of generators, we would have to have elements of arbitrarily large denominator. Further, we could fix an $N$ sufficiently large so that for a sequence of $n_i\ra\infty$, we would have some $r_i/n_i\in R(f)$ and $|r_i/n_i|\leq N$. (For example, given $r_1/n_1$, we can take $N=r_1/n_1$ by appropriately subtracting off multiples of $r_1/n_1$ from any other $r_i/n_i$.) But then we would have torsion-free representatives $\al^{r_i/n_i}$ satisfying
 $
  h(\bal^{r_i/n_i})\leq N\,h(\bal),
 $
 and as Lehmer irreducible representatives, each representative has the same degree $\delta(f)$, and thus we have an infinite number of algebraic numbers with bounded height and degree, contradicting Northcott's theorem.
\end{proof}
 The proof of Proposition \ref{prop:d-and-delta} will now proceed from the following series of lemmas:
\begin{lemma}
 Let $0\neq q\in\bQ$. Then
 $
 R(qf) = \frac{1}{q} R(f).
 $
\end{lemma}
\begin{proof}
 This is clear from the definition. 
\end{proof}
\begin{lemma}
 Let $f\in\cL$ with $R(f)=\bZ$ and let $p$ be a prime. Let $\beta$ be a torsion-free representative of $f$ and denote by $\beta^{1/p^n}$ any representative of the class of $\beta^{1/p^n}$ modulo torsion of minimal degree. Then
 \[
  \deg_{\bQ} \beta^{1/p^n} = p^n\,\deg_{\bQ} \beta=p^n \delta(\beta)\quad\text{for all}\quad n\in\bN.
 \]
\end{lemma}
\begin{proof}
 By choosing a representative $\beta$ in $\Qbar^\times$ of degree $\delta(f)$ we can say that $\deg_{\bQ} \beta^{1/p^n} \leq p^n\,\deg_{\bQ} \beta=\delta(f)$. Let us show that we cannot, in fact, do better if $R(f)=\bZ$. We proceed by induction. First observe that $\deg_{\bQ} \beta^{1/p}=p\deg_{\bQ} \beta$ because otherwise $1/p\in R(f)$, which contradicts our assumption. Suppose that $\deg_{\bQ} \beta^{1/p^{n-1}} =p^{n-1}\,\deg_{\bQ} \beta$, $\deg_{\bQ} \beta^{1/p^{n}} =p^{n}\deg_{\bQ} \beta$ but $\deg_{\bQ} \beta^{1/p^{n+1}} =p^{n}\,\deg_{\bQ} \beta$. Then we have the following tower of fields:
\[\xymatrix{
\bQ(\beta^{1/p^n})=\bQ(\beta^{1/p^{n+1}}) \ar@{-}^{p}[d]\\
\bQ(\beta^{1/p^{n-1}})} 
\]
Then over $\bQ(\beta^{1/p^{n-1}})$, $\beta^{1/p^{n+1}}$ is a root the polynomial in $\bZ[x]$ given by
\[
 G(x)=x^{p^2}-\beta^{1/p^{n-1}}=\prod_{i=1}^{p^2} (x-\zeta^i \beta^{1/p^{n+1}})
\]
where $\zeta$ denotes a primitive $p^2$th root of unity. But as $G(x)$ must have an irreducible factor $H(x)$ of degree $p$ over $\bQ(\beta^{1/p^{n-1}})$, the constant term of this polynomial is
\[
 H(0) = \zeta^m \beta^{1/p^n} \in \bQ(\beta^{1/p^{n-1}}),
\]
where $m\in\bZ$, and hence we have constructed a representative of $\beta^{1/p^n}$ that has degree $p^{n-1}\,\deg_{\bQ} \beta$, which contradicts the induction hypothesis. 
\end{proof}
\begin{lemma}\label{lemma:deg-beta-kn}
 Suppose $R(f)=\bZ$ and let $\beta$ be a torsion-free representative of $f$. Suppose we have an $n$th root of $\beta$, denoted $\beta^{1/n}$, which satisfies $\deg_{\bQ} \beta^{1/n}=n\deg_{\bQ} \beta$. Then $\deg_{\bQ} \beta^{k/n}=n \deg_{\bQ} \beta$ for all $k\in\bZ$ with $(k,n)=1$.
\end{lemma}
\begin{proof}
Suppose $\deg_{\bQ} \beta^{k/n}<\deg_{\bQ} \beta^{1/n}=n\,\deg_{\bQ} \beta$. Then $t=[\bQ(\beta^{k/n}):\bQ(\beta)]<n$ and $x^n-\beta^k$ has an irreducible factor of degree $t$ over $\bQ(\beta)$. But then by considering the constant term of this polynomial, we see that there is an $n$th root of unity $\zeta$ such that $\zeta\beta^{kt/n}\in\bQ(\beta)$ and hence has degree $\deg_{\bQ} \beta=\delta(f)$ and is Lehmer irreducible. But $kt/n\not\in\bZ$ since $(k,n)=1$ and $t$ is a proper divisor of $n$, thus $R(f)$ is strictly larger than $\bZ$, which contradicts our assumption. Therefore we must have $\deg_{\bQ} \beta^{k/n}=\deg_{\bQ} \beta^{1/n}$, and the proof is complete.
\end{proof}
\begin{lemma}
Suppose $f\in\cL$ with $R(\beta)=\bZ$ and let $\beta\in\Qbar^\times$ be a torsion-free representative. Suppose $n,m\in\bN$ are such that $(n,m)=1$ and that
\[
 \deg_{\bQ} \beta^{1/n} = n\,\deg_{\bQ} \beta\quad\text{and}\quad
 \deg_{\bQ} \beta^{1/m} = m\,\deg_{\bQ} \beta.
\]
Then $\deg_{\bQ} \beta^{1/mn}=mn\deg_{\bQ} \beta$ for an $nm$-th root of minimal degree $\beta^{1/mn}$.
\end{lemma}
\begin{proof}
 Choose representatives as in the proof of the lemma above. Choose $k,\ell\in\bZ$ such that $km+n\ell=1$ and thus $\beta^{k/n}\beta^{\ell/m}=\beta^{1/mn}$. By the above lemma, $\deg_{\bQ} \beta^{k/n}=\deg_{\bQ} \beta^{1/n}$ and likewise for $\beta^{\ell/m}$, and so since the degrees are relatively prime over $\bQ(\beta)$, we have the desired result:
\[\xymatrix@ur{
 \bQ(\beta^{k/n})\ar@{-}^m[r] & \bQ(\beta^{1/nm})\ar@{-}^n[d]\\
 \bQ(\beta)\ar@{-}^n[u]\ar@{-}_m[r] & \bQ(\beta^{\ell/m})}.\qedhere
\]
\end{proof}

Combining the above three lemmas, we now see that we have the proof of Proposition \ref{prop:d-and-delta}.

\subsection{Reduction to Lehmer irreducible numbers}
We will now show that we can reduce questions related to the $L^p$ Mahler measure to the set of Lehmer irreducible elements. We begin with two lemmas regarding the relationship between the projection operators $P_K$ and the degree functions $d$ and $\delta$ which will be used below:
\begin{lemma}\label{lemma:d-P-leq-d}
If $f\in\cF$ and $K\subset  K_f$, then $d(P_Kf)\leq d(f)$.
\end{lemma}
\begin{proof}
 Let $f=\fal$ and let $\al\in\Qbar^\times$ be a minimal degree representative of $f$, and choose $\ell\in\bN$ such that $\al^\ell$ is torsion-free. Then $\bQ(\al^\ell)= K_f$, so in particular, we see that
 \[
  K\subseteq  K_f\subseteq \bQ(\al).
 \]
 Observe that the norm $N_K^{K(\al)}$ from $K(\al)$ to $K$ is well-defined on the class $\bal\in\cG$. Since $(N_K^{K(\al)}\al)^{1/[K(\al):K]}$ is a representative of $({N_K^{K(\al)}\al})^{1/[K(\al):K]}$ modulo torsion, it follows from the fact that $N_K^{K(\al)}\al\in K$ that
 \begin{multline*}
  d(P_K f)\leq \deg_{\bQ}(N_K^{K(\al)}\al)^{1/[K(\al):K]}\leq[K(\al):K]\cdot[K:\bQ]\\
  =[\bQ(\al):\bQ]=d(f).\qedhere
 \end{multline*}
\end{proof}

\begin{lemma}\label{lemma:delta-P-leq-delta}
 If $K\in\cK$ and $K\subset K_f$ for $f\in\cF$, we have $\delta(P_Kf)\leq \delta(f)$.
\end{lemma}
\begin{proof}
 Since we can rescale $f$ without affecting either $\delta$ value, we can assume $f\in\cL$ so $d(f)=\delta(f)$. Let $F= K_f$. Then by Lemma \ref{lemma:d-P-leq-d} above, we have
 \[
  \delta(P_Kf)\leq d(P_Kf)\leq d(f)=\delta(f).\qedhere
 \]
\end{proof}
From the construction of $d$ above, it is easy to see that:
\begin{prop}
Let $m_p :\cF\ra [0,\infty)$ be given by $m_p(f) = d(f)\cdot \|f\|_p$. Fix $0\neq f\in\cF$. Then
\[
 m_p(f) = \min \{ (\deg_{\bQ}\al)\cdot h_p(\al) : \al\in\Qbar^\times, \fal=f\}.
\]
The right hand side of this equation is the minimum of the $L^p$ analogue of the usual logarithmic Mahler measure on $\Qbar^\times$ taken over all representatives of $f$ modulo torsion.
\end{prop}
We now prove the reduction to $\cL\subset\cF$:
\begin{prop}\label{prop:suffices-to-bound-on-L}
 Let $m_p(f) = d(f)\cdot \|f\|_p$. Then
 $
  m_p(\cF) = m_p(\cL),
 $
 so in particular, $\inf m_p(\cF\setminus\{0\})>0$ if and only if $\inf m_p(\cL\setminus\{0\})>0$.
\end{prop}
\begin{proof}
 Let $f\in\cF$ and $\ell=\ell(f)$. Then by Proposition \ref{prop:d-and-delta} we have $\delta(f)=d(\ell f)$ and $\ell\,\delta(f)=d(f)$, and thus
 \[
  m_p(\ell f) = \delta(f)\cdot \|\ell f\|_p = \ell\,\delta(f)\|f\|_p = d(f)\cdot\|f\|_p=m_p(f).\qedhere
 \]
\end{proof}

\begin{rmk}
 Proposition \ref{prop:suffices-to-bound-on-L}, which will be used below in the proof of Theorem \ref{thm:Lp-Lehmer-equiv}, is a key step in constructing equivalent statements of Lehmer's conjecture for heights which scale, such as $\delta\, h_p$ and particularly for the norms we will construct. Consider for example that if $\al = 2^{1/n}$ then $\delta(\fal)=1$ for all $n\in\bN$ and $h_1(2^{1/n}) = (2 \log 2)/n\ra 0$.
\end{rmk}

\subsection{Projection irreducibility}
In this section we introduce the last criterion which we will require to reduce the Lehmer conjectures to a small enough set of algebraic numbers to prove our main results.
\begin{defn}
We say $f\in\cF$ is \emph{projection irreducible} if $P_K(f) = 0$ for all proper subfields $K$ of the minimal field $K_f$. We denote the collection of projection irreducible elements by $\cP\subset\cF$.
\end{defn}

\begin{rmk}
 Notice that we cannot in general require that $P_K(f)=0$ for all $K\neq K_f$, as an element with a minimal field which is not Galois will typically have nontrivial projections to the conjugates of its minimal fields. See Remark \ref{rmk:why-not-decompose-K} above for more details.
\end{rmk}

We now prove that we can reduce questions about lower bounds on the Mahler measure $m_p$ to elements of $\cP$:
\begin{prop}\label{prop:reduction-to-T-irred}
 We have
 \[
  \inf_{f\in\cF\setminus\{0\}} m_p(f) > 0\quad \iff \quad \inf_{f\in\cP\setminus\{0\}} m_p(f)>0.
 \]
\end{prop}
\begin{proof}
 Let $f\in\cF$. Notice that for any $K\in\cK$ that by Lemma \ref{lemma:d-P-leq-d} we have 
 $d(\ovrln{P_K f})\leq d(f)$ and by Lemma \ref{prop:PK-continuous-wrt-Weil-norm} we have $h_p(\ovrln{P_K f})\leq h_p(f)$, so $m_p(\ovrln{P_K \al})\leq m_p(f)$. Let $\supp_{\cK}(f)=\{K\in\cK : P_K f\neq 0\}$. Notice that if $K\subset L$ and $K\in \supp_{\cK}(f)$, then $L\in \supp_{\cK}(f)$. Let $E$ denote the Galois closure of $ K_f$, and observe that $P_K f = P_K(P_E f) = P_{K\cap E} f$ by Lemma \ref{lemma:PKcommute}, so since we have only a finite number of subfields of $E$, we can write $\supp_{\cK}(f) = \bigcup_{i=1}^n [K_i,\ )$ where $ [K_i,\ ) = \{ L \in\cK : K_i\subseteq L\}$, and each $K_i\subseteq E$ is minimal in the sense that $[K_i,\  )\not\subseteq [K_j,\ )$ for all $i\neq j$. Thus, for each $i$, $P_{F} f = 0$ for all $F\subsetneqq K_i$, and so $P_{K_i}f\in \cP\setminus\{0\}$. Then $0<m_p(P_{K_i}f)\leq m_p(f)$, and so we have shown $\inf_{f\in\cP\setminus\{0\}} m_p(f)\leq \inf_{f\in\cF\setminus\{0\}} m_p(f)$. The reverse inequality is trivial.
\end{proof}

\section{The Mahler $p$-norm}\label{sect:mahler-p}

\subsection{The Mahler $p$-norms and proof of Theorem \ref{thm:Lp-Lehmer-equiv}}
We will now make use of our orthogonal decomposition \eqref{eqn:decomposition-by-degree} to define one of the main operators of our study. Let
\begin{equation}
\begin{split}
 M : \cF &\ra \cF\\
 f &\mapsto \sum_{n=1}^\infty n\, T^{(n)} f.
\end{split}
\end{equation}
The $M$ operator serves the purpose of allowing us to scale a function in $\cF$ by its appropriate degree while still being linear. As each element of $\cF$ has a finite expansion in terms of $T^{(n)}$ components, the above map is well-defined. Further, it is easily seen to be linear by the linearity of the $T^{(n)}$, and it is also a bijection. However, it is not a bounded operator with respect to any $L^p$ norm, as elements $f=T^{(n)} f$ can be found in the subspaces $T^{(n)}_K(\cF)$ for $K\in\cK^G$ of unbounded degree (otherwise by \eqref{eqn:decomposition-by-deg-and-field} all algebraic numbers would have finite Galois orbit modulo torsion, which is absurd), and since for such an element we have $Mf = n\,f=[K:\bQ]\cdot f$, we can conclude that the map $M$ is unbounded. In particular, $M$ is not well-defined $L^p(Y)$.

We define the Mahler $p$-norm on $\cF$ to be
\begin{equation}\label{eqn:mahler-norm-defn}
 \|f\|_{m,p} = \|Mf\|_p
\end{equation}
where $\|\cdot\|_p$ denotes the usual $L^p$ norm as defined above. Observe that this construction does in fact define a vector space norm, because $M$ is both linear and invertible as an operator taking $\cF\ra\cF$. We can complete $\cF$ with respect to $\|\cdot\|_{m,p}$ to obtain a real Banach space which we denote $\cF_{m,p}$. We are now ready to prove Theorem \ref{thm:Lp-Lehmer-equiv}, which we restate for the reader's convenience:
\begin{thm2} 
 For each $1\leq p\leq \infty$, equation \eqref{eqn:Lp-lehmer-conj} holds if and only if 
 \begin{equation}\tag{$**_p$}
  \|f\|_{m,p} \geq c_p > 0\quad\text{for all}\quad 0\neq f\in\cL\cap \cP\cap\cU
 \end{equation}
 where $\cL$ denotes the set of Lehmer irreducible elements, $\cP$ the set of projection irreducible elements, and $\cU$ the subspace of algebraic units. Further, for $1\leq p\leq q\leq\infty$, if equation \eqref{eqn:Lp-lehmer-conj} holds for $p$ then equation \eqref{eqn:Lp-lehmer-conj} holds for $q$ as well.
\end{thm2}
We recall the equation
\begin{equation}\tag{$*_p$} 
 m_p(\al)=(\deg_{\bQ}\al)\cdot h_p(\al)\geq c_p>0\quad\text{for all}\quad \al\in\Qbar^\times\setminus \Tor(\Qbar^\times).
\end{equation}
from Conjecture \ref{conj:lehmer-p} above is the $L^p$ analogue of the Lehmer conjecture.
\begin{proof}[Proof of Theorem \ref{thm:Lp-Lehmer-equiv}]
 First let us show that it suffices to bound $m_p(f)$ away from zero for $f\in \cL\cap\cP\cap\cU$. The reduction to $\cL$ was proven above in Proposition \ref{prop:suffices-to-bound-on-L}. Let $f\in\cF$. Let us first reduce to the set $\cU=\{f\in\cF : \supp_Y(f)\subseteq Y(\bQ,\infty)\}$. If $1\leq q<\infty$, observe that 
 \[
  h_q(f)=\|f\|_q = \bigg(\sum_{p\in M_{\bQ}} \| f|_{Y(\bQ,p)}\|_q^q \bigg)^{1/q}\geq \| f|_{Y(\bQ,p)}\|_q\geq \| f|_{Y(\bQ,p)}\|_1,
 \]
 since $Y(\bQ,p)$ is a space of measure $1$. Likewise, it is easy to see that 
 \[
 h_\infty(f)=\max_{p\in M_{\bQ}} \| f|_{Y(\bQ,p)}\|_\infty \geq \| f|_{Y(\bQ,p)}\|_\infty \geq \| f|_{Y(\bQ,p)}\|_1 
 \]
 for a specific rational prime $p$, so we can let $q=\infty$ as well. Let the rational prime $p$ be chosen above so that the $q$-norm is nonzero, which we can do if $f\not\in\cU$. Let $\al\in\Qbar^\times$ be a representative of minimal degree $d(f)$ for $f$. Then $\al$ has a nontrivial valuation over $p$, and since the product of $\al$ over all of its conjugates must be in $\bQ$, we know that we must have $\| f|_{Y(\bQ,p)}\|_1\geq (\log p) / d(f)$. Thus $h_q(f)\geq (\log 2)/d(f)$, so $m_q(f)\geq \log 2$ for $1\leq q\leq \infty$ if $f\not\in\cU$. Now it remains to show that we can reduce to the consideration of $\cP$ as well, but this now follows immediately from Proposition \ref{prop:reduction-to-T-irred} above.
 
 Now let $f\in \cL\cap\cP\cap\cU$, and we will show that $m_p(f)=\|f\|_{m,p}$, completing the proof of the equivalence. Observe that for such an element, by projection irreducibility, we must have $T^{(n)} f= f$ where $n=[ K_f:\bQ]$ and $K_f$ is the minimal field of $f$, as otherwise we could find a minimal $m<n$ such that $T^{(m)}f \neq 0$ and we could write $T^{(m)}f$ as a sum of elements belonging to subspaces $V_K$ for $K\subsetneqq  K_f$ of degree $[K:\bQ]\leq m$, and thus $f$ would have to have a nontrivial projection to a minimal such subfield, contradicting its projection irreducibility. Thus
 \[
  \|f\|_{m,p} = \|Mf\|_p = [ K_f:\bQ] \cdot \|f\|_p = \delta(f)h_p(f) = d(f)h_p(f) = m_p(f).
 \]
 where the second inequality follows from the fact that $f\in \cP$ and the fourth from the fact that $f\in\cL$. This completes the equivalence of the bounds. 
 
 To show that for $1\leq p\leq q\leq\infty$ the result for $p$ implies the result for $q$, we observe that having reduced the problem to the study of algebraic units $\cU$, that these numbers are of the form
 \[
  \cU = \{f\in\cF : \supp_Y(f)\subseteq Y(\bQ,\infty)\}
 \]
 and since $\lambda(Y(\bQ,\infty))=1$, we are reduced to the consideration of measurable functions on a probability space $(Y(\bQ,\infty),\lambda)$. But on such a space one has the usual inequality $\|f\|_p\leq \|f\|_q$ and thus $\|f\|_{m,p}=\|Mf\|_p\leq \|Mf\|_q=\|f\|_{m,q}$.
\end{proof}

Lastly, we note for its own interest:
\begin{prop}\label{prop:equivalent-to-classical-conjectures}
 Equation \eqref{eqn:Lp-lehmer-conj} for $p=1$ is equivalent to the Lehmer conjecture, and for $p=\infty$, \eqref{eqn:Lp-lehmer-conj} is equivalent to the Schinzel-Zassenhaus conjecture.
\end{prop}
\begin{proof}
  Since $h=2h_1$ it is obvious that $m_1=2m$ so we exactly have the statement of the Lehmer conjecture when $p=1$. Let us now show that when $p=\infty$, equation \eqref{eqn:Lp-lehmer-conj} is equivalent to the Schinzel-Zassenhaus conjecture. Recall that the house $\shorthouse\al = \max\{|\sigma\al| : \sigma:\bQ(\al)\hookrightarrow\bC\}$ where $|\cdot|$ denotes the usual Euclidean absolute on $\bC$. The Schinzel-Zassenhaus conjecture \cite{SchZas} states that for an algebraic integer $\al$, $(\deg_{\bQ}\al)\cdot \log \shorthouse{\al}$ is bounded away from zero by an absolute constant. Observe that by Smyth's well-known theorem \cite{Smy}, we have $m_1(\al)\geq c>0$ for an absolute constant $c$ if $\al$ is not reciprocal. Since $\|f\|_{m,\infty}\geq \|f\|_{m,1}=m_1(f)$ for the numbers under consideration, we see that if $\al$ is not reciprocal, then there is nothing more to show by the previous theorem. If $\al$ is reciprocal, then observe that $\al$ and $\al^{-1}$ are conjugate, and so $\house{\al}=\max\{\house{\al},\tallhouse{\al^{-1}}\}$, where $\max\{\house{\al},\tallhouse{\al^{-1}}\}$ is called the \emph{symmetric house}. Now, it is easy to see that $h_\infty(\al)=\log\max\{\house{\al},\tallhouse{\al^{-1}}\}$ is the logarithmic symmetric house of $\al$ for $\fal\in\cU$, so we do indeed recover the Schinzel-Zassenhaus conjecture when $p=\infty$.\footnote{We remark in passing that while $h_\infty$ agrees with the logarithmic symmetric house on $\cU$, $h_\infty$ seems to be a better choice for non-integers as well, as, for example, $h_\infty(3/2)=\log 3$ while the logarithmic symmetric house of $3/2$ is $\log(3/2)$.}
\end{proof}

\subsection{Explicit values}
We now evaluate the Mahler $p$-norms for two classes of algebraic numbers, surds and Salem numbers. Salem numbers are conjectured to be of minimal Mahler measure for the classical Lehmer conjecture. This is in part due to the fact that the minimal value for the Mahler measure known, dating back to Lehmer's original 1933 paper \cite{Lehmer}, is that of the Salem number called Lehmer's $\tau>1$, the larger positive real root of the irreducible polynomial $x^{10} + x^9 - x^7 - x^6 - x^5 - x^4 - x^3 + x + 1$. Here we show that, in fact, Salem numbers belong to the set $\cL\cap\cP\cap\cU$.

\subsubsection{Surds}
Recall that a \emph{surd} is a number $f\in\cF$ such that $\delta(f)=1$, which is equivalent to $[ K_f:\bQ]=1$ so $ K_f=\bQ$, and thus $f\in V_{\bQ}$. Now, $T_{\bQ}=P_{\bQ}$, and therefore, all surds are projection irreducible as they are fixed by $T_{\bQ}$. Thus, for $f$ with $f$ a surd,
\[
 \|f\|_{m,p} = \delta(f) \|f\|_p = \|f\|_p = h_p(f).
\]

\subsubsection{Pisot and Salem numbers}
We say that $f_\tau\in\cF$ is Pisot or Salem number if it has a representative $\tau\in\Qbar^\times$ which is a Pisot or Salem number, respectively. Recall that $\tau>1$ is said to be a Pisot number if $\tau$ is an algebraic integer whose conjugates in the complex plane all lie strictly within the unit circle, and that $\tau>1$ is a Salem number if $\tau$ is algebraic unit which is reciprocal and has all conjugates except $\tau$ and $\tau^{-1}$ on the unit circle in the complex plane (with at least one pair of conjugates on the circle).
\begin{prop}
 Every Pisot or Salem number $f_\tau$ is Lehmer irreducible, that is, $f_\tau\in\cL$.
\end{prop}
\begin{proof}
 It is easy to see that for a Pisot or Salem number $f_\tau$ and its given representative $\tau>1$, that $\house{\tau}=\tau$ and all other Galois conjugates $\tau'$ have $|\tau'|<|\tau|$. Therefore $\tau$ is Lehmer irreducible, since if $\delta(f_\tau)<\deg_{\bQ} \tau$, then each equivalence class modulo torsion would have more than one member, and in particular the real root $\tau>1$ would not uniquely possess the largest modulus, as $\zeta\tau$ would be a conjugate for some $1\neq \zeta\in \Tor(\Qbar^\times)$ which would have the same modulus, a contradiction. Since $f_\tau$ has a representative of degree $\delta(f_\tau)$, we have by definition $f_\tau\in\cL$.
\end{proof}
\begin{prop}
 Every Salem number $\tau$ is projection irreducible, that is, $f_\tau\in\cP$.
\end{prop}
\begin{proof}
 Suppose $f_\tau$ has its distinguished representative $\tau\in K^\times$, where $K=K_f=\bQ(\tau)$. Then there are precisely two real places of $K$, call them $v_1,v_2|\infty$, where $\tau$ has nontrivial valuation, and the remaining archimedean places are complex. By the definition of projection irreducibility, we need to show that $P_F(f_\tau) = 0$ for all $F\subsetneqq K$. Now, since $\lambda(Y(K,v_1))=\lambda(Y(K,v_2))=1/[K:\bQ]$, we know that for our subfield $F\subsetneqq K$, either $Y(K,v_1)\cup Y(K,v_2)\subseteq Y(F,w)$ for some place $w$ of $F$, in which case $P_F(f_\tau)=0$ because the two valuations sum to zero by the product formula, or else $v_1$ and $v_2$ lie over distinct places of $F$, call them $w_1$ and $w_2$. Then the algebraic norm $\beta=\N^K_F \tau$ has nontrivial valuations at precisely the two archimedean places $w_1,w_2$. Observe that $w_1,w_2$ must be real, as the completions are $\bQ_\infty=\bR\subset F_{w_i}\subset K_{v_i}=\bR$ for $i=1,2$. Thus $\beta$ must be a nontrivial Salem number or a quadratic unit. In either case, if we assume WLOG that $\log\|\beta\|_{w_1}>0$, observe that
 \[
  \beta = \|\beta\|_{w_1}
 \]
 But it is easy to see that
 \[
  \log \|\beta\|_{w_1} = \frac{1}{[K:F]}\log\|\tau\|_{v_1}
 \]
 and thus $\beta^{[K:F]} = \tau$. But this is a contradiction, as then the minimal field of $f_\beta$ must also be $K$, but $\beta\in F\subsetneqq K$.
\end{proof}

Thus, if $\tau>1$ is a Salem number, we have $f_\tau\in\cL\cap\cP$, so we can compute explicitly:
\begin{equation}
 \|f_\tau\|_{m,p} = \delta(f_\tau) \|f_\tau\|_p 
  = \delta(f_\tau)^{1-1/p} 2^{1/p} |\log \tau|.
\end{equation}
When $p=1$ this is, of course, twice the classical logarithmic Mahler measure of $\tau$, and when $p=\infty$, this is precisely the degree times the logarithmic house of $\tau$.

\subsection{The group $\Gamma$ and proof of Theorem \ref{thm:Gamma-equiv}}
We now expand the set which must be bounded away from 0 if the $L^p$ Lehmer conjecture is true to include the additive subgroup $\Gamma = \langle \cL\cap\cP\cap\cU\rangle$ and thus we establish Theorem \ref{thm:Gamma-equiv}.

\begin{lemma}\label{lemma:product-of-TILIs-in-same-field-is-TILI}
 Suppose $f,g\in\cL\cap\cP$ are projection irreducible for the same minimal field $K= K_f=K_g$. Then $f+g\in\cL\cap\cP$ as well, and if $f+g\neq 0$, then $K=K_{f+g}$ is the minimal field of $f+g$ as well.
\end{lemma}
\begin{proof}
 If $f+g=0$ then the problem is trivial as $0\in\cL\cap\cP$, so suppose $f+g\neq 0$ and let $K= K_f=K_g$. Then clearly $K_{f+g}\subseteq K$, and in fact, it easy to see that we must have equality, since if $F=K_{f+g}\subsetneqq K$ we would have $P_F(f+g)=f+g$, but $P_F(f+g) = P_F f + P_f g = 0 + 0 = 0$, a contradiction. Thus we have $f+g\in\cP$. Now choose torsion-free representatives $\al,\beta\in K^\times$ of $f,g$ respectively. It remains to show that $\al\beta$ has a torsion-free representative as well to show that the class $f+g$ is an element of $\cL$. For our chosen representatives, observe that the product is in $K$ as well. Let $\ell\in\bN$ be the minimal power to which we must raise $\al\beta$ to ensure it is torsion-free. If we can show that $\ell=1$, the proof will be complete. Observe that $K$ must also be the minimal field for the class of the torsion-free number $(\al\beta)^\ell$, and thus $K=\bQ((\al\beta)^\ell)$, so $(\al\beta)^\ell$ generates it and has full degree. But $\al\beta\in K$ and therefore generates it as well, and so $\al\beta$ is torsion-free, and the proof is complete.
\end{proof}

\begin{prop}
 The group $\Gamma$ is free abelian.
\end{prop}
\begin{proof}
 For each field $K\in\cK^G$, let $W_K = T_K(\cU)$. Then $W_K$ is a finite dimensional $\bQ$-vector space by Dirichlet's unit theorem. Each element $f$ of $\cP\cap\cU$, by definition of projection irreducible, belongs to $V_F\cap\cU$ for a unique minimal field $F\in\cK$. Then if $K\in\cK^G$ is the Galois closure of $F$, we have $f\in W_K$. Thus $\Gamma$ is generated by the elements of $(\cL\cap\cP\cap\cU)\cap W_K=\cL\cap\cP\cap W_K$ as $K$ ranges over $\cK^G$. In particular, observe that
 \begin{equation}
  \Gamma = \bigoplus_{K\in\cK^G} (\Gamma \cap W_K).
 \end{equation}
 Now, the set $\cL\cap \cP \cap W_K$ generates $\Gamma \cap W_K$. Further, $\cL\cap \cP \cap W_K$ can be viewed as a subset of the group of units modulo torsion of the field $K$, since the Lehmer irreducible representative of any $f\in W_K$ is well-defined in the multiplicative group of units of $K$ modulo torsion. But then the group $\langle \cL\cap \cP \cap W_K\rangle = \Gamma \cap W_K$ is free, since it is generated by a subset of a free abelian group of finite rank. Thus $\Gamma$ is a direct sum of finite rank free abelian groups and is free abelian itself.
\end{proof}

Let $\cU_{m,p}$ denote the completion of $\cU$ with respect to the Mahler $p$-norm $\|\cdot\|_{m,p}$. We now prove Theorem \ref{thm:Gamma-equiv}, which we recall here:
\begin{thm3}
 Equation \eqref{eqn:Lp-lehmer-conj} holds if and only if $\Gamma\subset \cU_{m,p}$ is closed.
\end{thm3}
\begin{proof}
 By Proposition \ref{prop:suffices-to-bound-on-L} and the argument of Theorem \ref{thm:Lp-Lehmer-equiv} regarding reducing to units, we know that \eqref{eqn:Lp-lehmer-conj} holds if and only if there exists a constant $c_p$ such that $m_p(f)\geq c_p>0$ for all $f\in\cL\cap\cU$. By the fact that $P_K$ is a norm one projection with respect to the $L^p$ norm (Propositon \ref{prop:PK-continuous-wrt-Weil-norm}), and the fact that it commutes with the $T^{(n)}$ operators (Proposition \ref{prop:PnPK-commute}), we see that it commutes with the $M$ operator well, and therefore, by the definition of the Mahler norm,
 \[
  \|P_K f\|_{m,p} = \|M P_K f\|_{p}=\|P_K(Mf)\|_{p}\leq \|Mf\|_p = \|f\|_{m,p}.
 \]
 Now, as a free abelian additive subgroup of the separable Banach space $\cU_{m,p}$, $\Gamma$ is discrete if and only it is closed. If $\Gamma$ is discrete, then since $\cL\cap\cP\cap\cU\subset \Gamma$ we have the desired result by Theorem \ref{thm:Lp-Lehmer-equiv}. Suppose on the other hand that we know that equation \eqref{eqn:Lp-lehmer-conj} holds. Let $f\in\Gamma$. Then $f = \sum_{i=1}^n g_i$ where the $g_i\in \cL\cap\cP\cap\cU$. As projection irreducible elements, each $g_i$ has a unique minimal field $K_{g_i}$ associated to it and has no nontrivial projections to any proper subfields of $K_{g_i}$. The set of fields $A=\{K_{g_i} : 1\leq i\leq n\}\subset\cK$ is finite, and therefore, there must exist an element $K\in A$ which is minimal in this set, that is, there is no element $F\in A$ such that $F\subsetneqq K$. We can assume $P_K(f)\neq 0$, otherwise, the elements $g_i$ with $K_{g_i}=K$ would sum to zero and we could remove them from the sum expressing $f$ without changing the value. By the above inequality for $P_K$, we have $\|P_K f\|_{m,p}\leq \|f\|_{m,p}$. But 
 \[
  P_K f = \sum_{\substack{1\leq i\leq n\\ K_{g_i}=K}} g_i,
 \]
 and therefore by Lemma \ref{lemma:product-of-TILIs-in-same-field-is-TILI} above, $P_Kf \in \cL\cap\cP\cap\cU$. Then by assumption and by Theorem \ref{thm:Lp-Lehmer-equiv}, we have an absolute constant $c_p$ such that
 \[
  \|f\|_{m,p} \geq \|P_K f\|_{m,p}\geq c_p > 0.
 \]
 Thus $\Gamma$ is indeed discrete and therefore closed, as claimed.
\end{proof}

\subsection{The Mahler $2$-norm and proof of Theorem \ref{thm:mahler-2-norm}}
Recall that we define the Mahler $2$-norm for $f\in\cF$ to be:
\[
 \| f \|_{m,2} = \| Tf\|_2 = \bigg\| \sum_{n=1}^\infty n\,T^{(n)} f \bigg\|_2.
\]
The goal of this section is to prove Theorem \ref{thm:mahler-2-norm}, which we recall here for the convenience of the reader:
\begin{thm4}
 The Mahler $2$-norm satisfies
 \[
  \|f\|_{m,2}^2 = \sum_{n=1}^\infty n^2\,\|T^{(n)}(f)\|_2^2= \sum_{K\in\cK^G}\sum_{n=1}^\infty n^2\,\|T^{(n)}_K(f)\|_2^2.
 \]
 Further, the Mahler $2$-norm arises from the inner product
 \[
  \langle f,g\rangle_m = \langle Mf,Mg\rangle = \sum_{n=1}^\infty n^2\, \langle T^{(n)} f,T^{(n)}g\rangle
  = \sum_{K\in\cK^G} \sum_{n=1}^\infty n^2\, \langle T_K^{(n)} f,T_K^{(n)}g\rangle
 \]
where $\langle f,g\rangle =\int_Y fg\,d\lambda$ denotes the usual inner product in $L^2(Y)$, and therefore the completion $\cF_{m,2}$ of $\cF$ with respect to the Mahler $2$-norm is a Hilbert space.
\end{thm4}
\begin{proof}[Proof of Theorem \ref{thm:mahler-2-norm}]
 The first part of the theorem follows easily from the fact that the $T^{(n)}_K$ form an orthogonal decomposition of $\cF$. Indeed, for $f\in\cF$, we have:
 \[
  \|f\|_{m,2}^2 = \bigg\|\sum_{K\in\cK^G}\sum_{n=1}^\infty n\,T^{(n)}_K(f)\bigg\|_2^2
  = \sum_{K\in\cK^G}\sum_{n=1}^\infty n^2\,\|T^{(n)}_K(f)\|_2^2.
 \]
 The above sums are, of course, finite for each $f\in\cF$. That the specified inner product $\langle f,g\rangle _m$ defines this norm is then likewise immediate. Therefore, the completion of $\cF$ with respect to the norm $\|\cdot\|_{m,2}$ is a Hilbert space, as claimed.
\end{proof}

Lastly, we note that
$
\|\cdot\|_{m,2} \leq \delta h_2\leq m_2.
$
The authors suspect that this inequality is not true for general $p\neq 2$, but we know of no examples proving such a result. To see that the desired inequality holds for $p=2$, let us recall (Proposition \ref{prop:Tn-expansion}) that for a given $f\in\cF$, we have an expansion into degree $n$ components given by $f = T^{(1)} f + \cdots + T^{(N)} f$ with $T^{(N)} f\neq 0$. Then observe that $\delta(f) \geq N$, for otherwise, $T^{(N)} f=0$ since $f$ itself would have $[ K_f:\bQ]=n<N$ and thus 
 $f\in V^{(n)}$, and so it would have no essential projection to $V^{(N)}$. Thus
\begin{multline*}
 \|f\|_{m,2} = \left(\sum_{n=1}^N n^2 \|T^{(n)} f\|_2^2\right)^{1/2}
 \leq \left(\sum_{n=1}^N N^2 \|T^{(n)} f\|_2^2\right)^{1/2}\\ = N \|f\|_2 \leq \delta(f)\|f\|_2=\delta h_2(f).
\end{multline*}
That $\delta h_2\leq m_2=d\,h_2$ follows from the inequality $\delta\leq d$.

\end{document}